\documentclass[a4paper,10pt,reqno]{amsart}
\usepackage{graphicx}
\usepackage{geometry}
\usepackage{float}
\usepackage{verbatim}
\usepackage{setspace}
\usepackage{siunitx}
\usepackage{paralist}
\usepackage{enumitem}
\usepackage{xcolor}
\usepackage{amsfonts,amsmath,amssymb,bbm,amsthm,amsbsy}
\usepackage[utf8]{inputenc}
\usepackage{lipsum}
\usepackage[english]{babel}
\usepackage{booktabs}
\usepackage{array}
\usepackage{subfig}
\usepackage{caption}
\usepackage{mathtools}
\usepackage{hyperref}

\usepackage{stmaryrd}

\newenvironment{frcseries}{\fontfamily{frc}\selectfont}{}

\makeatletter
\newtheorem*{rep@theorem}{\rep@title}
\newcommand{\newreptheorem}[2]{
\newenvironment{rep#1}[1]{
 \def\rep@title{#2 \ref{##1}}
 \begin{rep@theorem}}
 {\end{rep@theorem}}}
\makeatother

\newtheorem{thm}{Theorem}[section]
\newtheorem{lemma}[thm]{Lemma}
\newtheorem{prop}[thm]{Proposition}
\newtheorem{corr}[thm]{Corollary}

\newtheorem*{thm*}{Theorem}
\newtheorem*{lemma*}{Lemma}
\newtheorem*{prop*}{Proposition}
\newtheorem*{corr*}{Corrolary}
\newtheorem*{claim*}{Claim}

\theoremstyle{remark}
\newtheorem{rmk}[thm]{Remark}

\newtheorem*{rmk*}{Remark}
\newtheorem*{conj*}{Conjecture}
\newtheorem*{quest*}{Question}

\theoremstyle{definition}
\newtheorem{defn}[thm]{Definition}
\newtheorem{exmp}[thm]{Example}

\newtheorem*{defn*}{Definition}
\newtheorem*{exmp*}{Example}

\newreptheorem{theorem}{Theorem}
\newreptheorem{corollary}{Corollary}
\newreptheorem{proposition}{Proposition}

\usepackage{fancyhdr}
\usepackage{tikz}
\usepackage{tikz-cd}

\usepackage{IEEEtrantools}

\newenvironment{equ*}[1]{\begin{IEEEeqnarray*}{#1}}{\end{IEEEeqnarray*}}

\newcommand{\F}{\mathbb{F}}

\newcommand{\Z}{\mathbb{Z}}
\newcommand{\N}{\mathbb{N}}

\newcommand{\Nc}{\mathcal{N}}

\newcommand{\A}{\mathcal{A}}

\newcommand{\C}{\mathcal{C}}

\newcommand{\T}{\mathcal{T}}
\newcommand{\K}{\mathcal{K}}

\DeclareFontFamily{U}{mathx}{}
\DeclareFontShape{U}{mathx}{m}{n}{<-> mathx10}{}
\DeclareSymbolFont{mathx}{U}{mathx}{m}{n}
\DeclareMathAccent{\widecheck}{0}{mathx}{"71}

\newcommand{\inj}{\hookrightarrow}
\newcommand{\sur}{\twoheadrightarrow}

\newcommand{\kG}{k\llbracket G\rrbracket}

\DeclareMathOperator{\Del}{\Delta}
\DeclareMathOperator{\cotimes}{\,\widehat{\otimes}\,}
\DeclareMathOperator{\crotimes}{\,\widehat{\otimes}_\mathit{R}\,}

\DeclareMathOperator{\botimes}{\otimes^\blacksquare}
\DeclareMathOperator{\brotimes}{\otimes^\blacksquare_\mathcal{R}}

\DeclareMathOperator{\Hom}{Hom}
\DeclareMathOperator{\Bil}{Bil}

\DeclareMathOperator{\UHom}{\underline{\textup{Hom}}}
\DeclareMathOperator{\RHom}{\mathit{R}\textup{Hom}}
\DeclareMathOperator{\RUHom}{\mathit{R}\underline{\textup{Hom}}}

\DeclareMathOperator{\Ext}{Ext}
\DeclareMathOperator{\UExt}{\underline{\textup{Ext}}}
\DeclareMathOperator{\Tor}{Tor}

\DeclareMathOperator{\colim}{colim}

\newcommand{\Set}{\mathbf{Set}}
\newcommand{\Grp}{\mathbf{Grp}}
\newcommand{\Ab}{\mathbf{Ab}}
\newcommand{\Ring}{\mathbf{Ring}}
\newcommand{\ModR}{\mathbf{Mod}(R)}

\newcommand{\TAb}{\mathbf{TAb}}

\newcommand{\TModR}{\mathbf{TMod}(R)}
\newcommand{\Pro}{\mathbf{Pro}}
\newcommand{\CHED}{\mathbf{CHED}}

\newcommand{\PAb}{\mathbf{PAb}}
\newcommand{\PRing}{\mathbf{PRing}}
\newcommand{\PModR}{\mathbf{PMod}(R)}
\newcommand{\DModR}{\mathbf{DMod}(R)}
\newcommand{\CSet}{\mathbf{CondSet}}

\newcommand{\CAb}{\mathbf{CondAb}}
\newcommand{\CRing}{\mathbf{CondRing}}
\newcommand{\CModR}{\mathbf{CondMod}(R)}
\newcommand{\CModRR}{\mathbf{CMod}(\mathcal{R})}
\newcommand{\SAb}{\mathbf{SolidAb}}
\newcommand{\SR}{\mathbf{Solid}(\mathcal{R})}

\newcommand{\R}{\mathcal{R}}
\newcommand{\M}{\mathcal{M}}
\newcommand{\G}{\mathcal{G}}
\newcommand{\UR}{\underline{R}}
\newcommand{\UM}{\underline{M}}

\title{Profinite and Solid Cohomology}
\author{Jiacheng Tang}
\thanks{Email: jiacheng.tang@postgrad.manchester.ac.uk, University of Manchester}

\calclayout
\begin{document}
\maketitle

\begin{abstract}
Solid abelian groups, as introduced by Dustin Clausen and Peter Scholze, form a subcategory of all condensed abelian groups satisfying some ``completeness" conditions and having favourable categorical properties. Given a profinite ring $R$, there is an associated condensed ring $\underline{R}$ which is solid. We show that the natural embedding of profinite $R$-modules into solid $\underline{R}$-modules preserves $\mathrm{Ext}$ and tensor products, as well as the fact that profinite rings are analytic.
\end{abstract}

\section{Introduction}
\label{intro}

Condensed mathematics was recently introduced by Dustin Clausen and Peter Scholze as a framework to combine algebra and topology (see \cite{condensed}). A first observation is that topological abelian groups do not form an abelian category, but condensed abelian groups (and condensed modules) do. In fact, the category of condensed abelian groups (or condensed modules) satisfies the same Grothendieck's (AB) axioms as the category $\Ab$ of abelian groups.

Profinite groups and modules form an important class of topological groups and modules. Given a profinite ring $R$, the category $\PModR$ of profinite $R$-modules is abelian, but does not in general have exact coproducts. Thus, a natural question to ask is whether we can replace profinite modules with some condensed modules with better categorical properties.

The articles \cite{condensed}, \cite{analytic} and \cite{solidcoh} contain some of our results specialised to the case when the profinite ring is $\Z/p^n\Z$ or $\Z_p$, and we will prove that those statements hold for all profinite rings. Moreover, \cite{solidcoh} contains results about solid group cohomology (as its title hints at), whilst we shall deal with the more general case of solid ring cohomology, which will allow us to recover certain results in \cite{solidcoh}. The main aims of this paper are to show that condensed mathematics (specifically the solid theory) provides a good framework for studying profinite modules (see Theorem \ref{mainthm}), and to set up the necessary foundations for further research in the area.

Let $\CHED$ denote the category of compact Hausdorff extremally disconnected spaces (where a space is \emph{extremally disconnected} if the closure of every open set is open). Recall from \cite{condensed} that a \emph{condensed set/group/ring/\ldots} is a functor $$\mathcal{T}\colon\CHED^{\mathrm{op}}\to\Set/\Grp/\Ring/\ldots$$ such that $\mathcal{T}(\varnothing)=*$ (the terminal object) and such that for any $S_1, S_2\in\CHED$, the natural map $\mathcal{T}(S_1\sqcup S_2)\to\mathcal{T}(S_1)\times\mathcal{T}(S_2)$ is a bijection. Let $\CAb$ denote the category of condensed abelian groups. For $\R$ a condensed ring, we define a \emph{condensed $\R$-module} to be a condensed abelian group $\M$ with a natural transformation $\R\times\M\Rightarrow\M$ making $\M(S)$ a usual $\R(S)$-module for each $S\in\CHED$ (Definition \ref{condmoda}), or equivalently it is an $\R$-module object in $\CAb$ with the standard tensor product (Proposition \ref{condmod2a}). The category $\CModRR$ of condensed $\R$-modules has all small limits and colimits, is abelian satisfying the same Grothendieck's (AB) axioms as $\Ab$, and is generated by compact projective objects (Theorems \ref{thm3a} and \ref{thm4a}).

We establish some of these elementary results in the appendix, which can be safely skipped by those who are comfortable with basic condensed mathematics and sheaf theory. This part is by no means original and there should be plenty of literature on the subject. The appendix is mainly present as a reference for newcomers (which the author would have loved to read when he was first introduced to condensed mathematics a while ago). For condensed mathematics one may obviously refer to \cite{condensed}, or to \cite{land} for a gentler introduction, while for the theory of sheaves one may refer to \cite{macsheaves}. We will assume knowledge of basic category theory (refer to \cite{maccat}) and homological algebra (refer to \cite{weibel}).

Let $R$ be a profinite ring. The category $\PModR$ of profinite modules over $R$ is abelian with good limits (e.g.\ it is (AB3*), (AB4*), (AB5*)) and has enough projectives, but in general the coproduct is not exact (though it always exists). There is a duality (i.e.\ contravariant equivalence) between the category $\PModR$ of profinite $R$-modules and the category $\DModR$ of discrete (topological) $R$-modules, via the \emph{Pontryagin Duality functor} $\Hom(-,\mathbb{R}/\Z)$. This means that in practice, given any theorem for profinite modules, there is a dual one for discrete modules. In particular, $\DModR$ is an abelian category satisfying (AB3), (AB4), (AB5) and has enough injectives, but in general the product is not exact. For an introduction to profinite groups and modules, one may refer to \cite{profinite}.

There is a condensed ring $\underline{R}$ naturally associated to $R$, defined via $\underline{R}(S)=C(S,R)=\{\text{continuous maps }S\to R\}$ for $S\in\CHED$. We then have two canonical additive functors $\PModR\to\mathbf{CMod}(\underline{R})$ and $\DModR\to\mathbf{CMod}(\underline{R})$, both defined via $M\mapsto\underline{M}=C(-,M)$, which are fully faithful and exact (Theorem \ref{PDfullsub}), so we can view $\PModR$ and $\DModR$ as full abelian subcategories of $\mathbf{CMod}(\underline{R})$.

The above embeddings $C^P\colon\PModR\inj\mathbf{CMod}(\underline{R})$ and $C^D\colon\DModR\inj\mathbf{CMod}(\underline{R})$, which we call the \emph{condensation functors}, are not quite sufficient, as we will need \emph{solidity} (Definition \ref{solid}). For $\R$ a condensed ring, the category $\SR$ of solid $\R$-modules is an abelian subcategory of $\CModRR$ closed under all limits, colimits and extensions, and also generated by compact projective objects (Corollary \ref{solidcor}). Moreover, the essential images of the condensation functors $C^P$ and $C^D$ lie in $\mathbf{Solid}(\underline{R})$ (Lemma \ref{brinklemma1}). The importance of this is that the embedding $\PModR\inj\mathbf{Solid}(\underline{R})$ preserves $\Ext$ and tensor products (note that we do have a natural tensor product $-\otimes^\blacksquare_{\underline{R}}-=(-\otimes_{\underline{R}}-)^\blacksquare$ on $\mathbf{Solid}(\underline{R})$; see page \pageref{solidt}), which is our main theorem:

\begin{thm}[Theorem \ref{thm5}(ii), (iii) and Proposition \ref{rtensors}(iii)]\label{mainthm}
Let $R$ be a profinite ring and $M$ be a profinite $R$-module.
\begin{enumerate}[label=(\roman*)]
\item If $N$ is a profinite $R$-module, then $$\Ext^*_{\mathbf{Solid}(\underline{R})}(\underline{M},\underline{N})=\Ext^*_R(M,N).$$
\item If $N$ is a profinite $R$-module, then $$\underline{M}\otimes^\blacksquare_{\underline{R}}\underline{N}=\underline{M\crotimes N}.$$
\item If $N$ is a profinite $R$-module and $M$ is of type $FP_\infty$ (i.e.\ $M$ has a projective resolution in $\PModR$ with each term topologically finitely generated), or if $N$ is a discrete $R$-module, then the condensed version of (i) holds i.e.\ $$\UExt^*_{\mathbf{Solid}(\underline{R})}(\underline{M},\underline{N})=\underline{\Ext^*_R(M,N)},$$ where the $\Ext^*_R(M,N)$ have their natural topologies (see Theorem \ref{thm5} for the precise statement).
\end{enumerate}
\end{thm}

In (i), $\Ext^*_R$ means $\Ext^*_{\PModR}$, and either of the $\Ext^*$ can be defined using Yoneda extensions, since $\mathbf{Solid}(\underline{R})$ and $\PModR$ are abelian categories. Both categories have enough projectives, so we can also calculate $\Ext^*$ by projectively resolving the first variable. On the other hand, for $M$ profinite and $N$ discrete in (iii), $\Ext^*_R(M,N)$ is calculated by projectively resolving $M$ in profinite modules or injectively resolving $N$ in discrete modules or both. The $\UExt^*$ appearing on the left hand side of (iii) is the condensed version of $\Ext^*$: see the paragraph before Example \ref{uhomeg} for the definition of the condensed Hom $\UHom$, and see the paragraph after Lemma \ref{solidhoml} for the definition of $\UExt^*$.

Despite the condensation functors $C^P$, $C^D$ and the category $\mathbf{Solid}(\underline{R})$ having favourable properties (e.g.\ $\mathbf{Solid}(\underline{R})$ has exact coproducts though $\PModR$ does not), it is unclear to the author how much we can use the condensed category to prove results about profinite modules themselves.

In Section \ref{sec1}, we will fix notations and state some of the elementary results from the appendix (without proof). In Section \ref{sec2}, we will recall the notion of solid abelian groups and modules and prove the above theorem. The key proposition is Proposition \ref{limtensor} which roughly states that solid tensor products commute with products when profinite modules are involved. The section ends with miscellaneous items such as a quick discussion of condensed group rings (Proposition \ref{gmodsame}). To make Section \ref{sec2} more accessible, we shall postpone all discussions of derived categories to the next section. In Section \ref{sec3}, we recall the definition of analytic rings from \cite{condensed} and prove that profinite rings are analytic, which is an easy consequence of the key proposition. Indeed, \cite[Proposition 7.9]{condensed} shows that $\Z_p$ (with its natural free solid modules) is analytic, and its proof easily extends to all profinite rings.

Remark: whenever we write ``=" in this paper, we mean isomorphic, usually canonically isomorphic (or equivalent in the case of categories), e.g.\ instead of saying that the natural map $\mathcal{T}(S_1\sqcup S_2)\to\mathcal{T}(S_1)\times\mathcal{T}(S_2)$ is a bijection in the definition of a condensed set, we might write $\mathcal{T}(S_1\sqcup S_2)=\mathcal{T}(S_1)\times\mathcal{T}(S_2)$.

Convention: all rings have a 1 but are not necessarily commutative.

\subsection*{Acknowledgements}
The author would like to thank his supervisor Peter Symonds for his constant guidance, the following people (in alphabetical order) for helpful discussions on the subject matter: Matthew Antrobus, Rudradip Biswas, Emma Brink, Dustin Clausen, Peter Kropholler, Sof\'ia Marlasca Aparicio, Juan Esteban Rodr\'iguez Camargo, and Thomas Wasserman, as well as the following people (in alphabetical order) for reading drafts of this paper and giving solid feedback: Matthew Antrobus, Anton Farmar, Daniel Heath, and Gregory Kendall. The author would like give special thanks to Peter Scholze for his consistently swift email replies and condensed yet effective answers.

The author would also like to thank the anonymous referees for their useful comments. This work was supported by the EPSRC Doctoral Training Partnership [grant number \linebreak EP/W524347/1].

\section{Preliminaries}
\label{sec1}

The aim of this section is mostly to fix notations and remind the reader of basic results regarding condensed modules. See the appendix for proofs and more details.

Let $\mathbf{TGrp}$, $\mathbf{TAb}$, $\mathbf{TRing}$ and $\mathbf{TMod}(R)$ denote respectively the categories of topological groups, topological abelian groups, topological rings, and topological modules over a topological ring $R$. Further, let $\mathbf{Pro}$, $\mathbf{PGrp}$, $\mathbf{PAb}$, $\mathbf{PRing}$ and $\mathbf{PMod}(R)$ denote respectively the categories of profinite spaces, profinite groups, profinite abelian groups, profinite rings, and profinite modules over a topological ring $R$. We also write $\CHED$ for the category of compact Hausdorff extremally disconnected spaces. Recall the following from \cite{condensed}:

\begin{defn}\label{condset}
A \emph{condensed set/group/ring/\ldots} is a sheaf of sets/groups/rings/\ldots on the site of extremally disconnected spaces, with finite jointly surjective families of maps as covers. Equivalently, a condensed set/group/ring/\ldots is a functor $$\mathcal{T}\colon\CHED^{\mathrm{op}}\to\Set/\Grp/\Ring/\ldots$$ such that $\mathcal{T}(\varnothing)=*$ and for any $S_1, S_2\in\CHED$, the natural map $\mathcal{T}(S_1\sqcup S_2)\to\mathcal{T}(S_1)\times\mathcal{T}(S_2)$ is a bijection.

Given a condensed set/group/ring/\ldots $\mathcal{T}$, we call the collection of global sections $\mathcal{T}(*)$ its \emph{underlying set/group/ring/\ldots}.
\end{defn}

Let $\mathbf{CondSet}$, $\mathbf{CondGrp}$, $\mathbf{CondAb}$, $\mathbf{CondRing}$ and $\mathbf{CondMod}(R)$ denote respectively the (large) categories of condensed sets, condensed groups, condensed abelian groups, condensed rings, and condensed modules over an abstract ring $R$ i.e.\ the codomain of $\mathcal{T}$ is $\ModR$. Given a condensed ring $\R$, let $\CModRR$ denote the category of condensed $\R$-modules. (Note that we're not using the possibly better notation $\mathbf{CondMod}(\mathcal{R})$ to avoid confusion with the category $\CModR$ of condensed modules over an abstract ring $R$.)

Given $T$ a topological space/group/ring/\ldots, there is a natural way to associate to it a condensed set/group/ring/\ldots (\cite[Example 1.5]{condensed}), namely we define the functor $\underline{T}$ via $\underline{T}(S)=C(S,T)$, the set of continuous maps from $S\in\CHED$ to $T$. The group/ring/\ldots structure on $C(S,T)$ is pointwise induced by that of $T$. This indeed makes $\underline{T}$ a sheaf. Now let $R$ be a topological ring and $M$ a topological module over $R$, so we have $\underline{R}\in\CRing$ and $\underline{M}\in\CModR$, where the last $R$ is viewed as an abstract ring. Then $\underline{M}$ is also naturally a condensed $\underline{R}$-module.

\begin{thm}[Theorem \ref{thm1a}]\label{thm1}
Let $R$ be an abstract ring. The category $\CModR$ of ($\kappa$-small) condensed modules over $R$ has all small limits and colimits and is an abelian category satisfying (AB3), (AB4), (AB5), (AB6), (AB3*) and (AB4*).
\begin{proof}
Exactly the same as in \cite[Theorem 2.2]{condensed}. The main point is that limits and colimits are computed pointwise here (limits are generally computed pointwise for abelian sheaves on a site, but colimits usually need to be further sheafified).
\end{proof}
\end{thm}

\begin{thm}[Theorem \ref{thm2a}]\label{thm2}
Let $R$ be an abstract ring. The category $\CModR$ of ($\kappa$-small) condensed modules over $R$ is generated by compact projective objects, that is, for any $\M\in\CModR$ there are some compact projective objects $\mathcal{P}_i\in\CModR$ and a surjection $\bigoplus\mathcal{P}_i\sur\M$. In particular, $\CModR$ has enough projectives.
\begin{proof}
Exactly the same as in \cite[Theorem 2.2]{condensed}. Note that the forgetful functor $\CModR\to\CSet$ has a left adjoint $\mathcal{T}\mapsto R[\mathcal{T}]$, where $R[\mathcal{T}]$ is the sheafification of the functor that sends $S\in\CHED$ to the abstract free $R$-module $R[\mathcal{T}(S)]$. The compact projective generators are then given by $R[\underline{S}]$ for $S\in\CHED$.
\end{proof}
\end{thm}

The category $\CAb$ of condensed abelian groups has a tensor product $\otimes$ (see \cite[page 13]{condensed}). Specifically, given $\M,\mathcal{N}\in\CAb$, $\M\otimes\mathcal{N}$ is defined as the sheafification of $S\mapsto\M(S)\otimes_{\Z}\mathcal{N}(S)$. Just like for abstract abelian groups, giving a morphism $\M\otimes\mathcal{N}\Rightarrow\mathcal{K}$ in $\CAb$ is equivalent to giving a bilinear map $\M\times\mathcal{N}\Rightarrow\mathcal{K}$ i.e.\ a natural transformation $\M\times\mathcal{N}\Rightarrow\mathcal{K}$ in $\CSet$ such that for each $S\in\CHED$, $\M(S)\times\mathcal{N}(S)\to\mathcal{K}(S)$ is bilinear in $\Ab$. The tensor product $\otimes$ makes $\CAb$ a monoidal category with unit $\underline{\Z}$, the constant sheaf with value $\Z$.

\begin{thm}[Theorem \ref{thm3a}]\label{thm3}
Let $\R$ be a condensed ring. The category $\CModRR$ of ($\kappa$-small) condensed $\R$-modules has all small limits and colimits and is an abelian category satisfying (AB3), (AB4), (AB5), (AB6), (AB3*) and (AB4*).
\end{thm}

\begin{thm}[Theorem \ref{thm4a}]\label{thm4}
Let $\R$ be a condensed ring. The category $\CModRR$ of ($\kappa$-small) condensed $\R$-modules is generated by compact projective objects. In particular, $\CModRR$ has enough projectives.
\end{thm}

Remark: again, all limits and colimits in $\CModRR$ are computed pointwise. The compact projective generators are $\R\otimes\Z[\underline{S}]$ for $S\in\CHED$.

Given a condensed ring $\R$, a right $\R$-module $\M$ and a left $\R$-module $\mathcal{N}$, there is a natural way to define their tensor product over $\R$, namely as the presheaf of abelian groups $(\M\otimes_{\R}\mathcal{N})(S)=\M(S)\otimes_{\R(S)}\mathcal{N}(S)$ for $S\in\CHED$. Note that this is actually a sheaf, so the pointwise tensor product $\M\otimes_{\R}\mathcal{N}$ is a condensed abelian group.

As explained on \cite[page 13]{condensed}, given $\M,\mathcal{N}\in\CAb$, the abelian group of morphisms $\Hom(\M,\mathcal{N})$ can be enriched to a condensed abelian group $\UHom(\M,\mathcal{N})$ by defining, for $S\in\CHED$, $\UHom(\M,\mathcal{N})(S)=\Hom_\CAb(\Z[\underline{S}]\otimes\M,\mathcal{N})\in\Ab$. This defines an internal Hom satisfying the Hom-tensor adjunction $\Hom_\CAb(\mathcal{K},\UHom(\M,\mathcal{N}))=\Hom_\CAb(\mathcal{K}\otimes\M,\mathcal{N})$ (for a proof, see \cite[Proposition 5.5]{land}). The same argument works for general condensed rings, which is explained in the appendix. For $\R$ a condensed ring, we will write $\Hom_\R$ to mean $\Hom_{\CModRR}$, similar to the abstract case.

\begin{prop}[Proposition \ref{homtensora}]\label{homtensor}
Let $\R,\R'$ be condensed rings, $\K$ a right $\R$-module, $\M$ an $\R$-$\R'$-bimodule, and $\mathcal{N}$ a right $\R'$-module. Then we have the Hom-tensor adjunction $$\Hom_{\R}(\mathcal{K},\UHom_{\R'}(\M,\mathcal{N}))=\Hom_{\R'}(\mathcal{K}\otimes_\R\M,\mathcal{N}).$$
\end{prop}

\section{Solid Modules and (Co)homology}
\label{sec2}

\subsection{Embedding Topological Modules into Condensed Modules}
\label{sub1}

Let $R$ be a topological ring. Then the category $\TModR$ of topological $R$-modules may not be abelian. For example, the category $\TAb$ of topological abelian groups is not abelian. This is an obstacle to using homological algebra. However, recall from Lemma \ref{topmoda} that there is an additive functor $C\colon\TModR\to\mathbf{CMod}(\underline{R})$ sending a topological $R$-module $M$ to the condensed $\underline{R}$-module $\underline{M}=C(-,M)$. We call $C$ the \emph{condensation functor} and will sometimes use $C$ instead of $\underline{(-)}$ or $M\mapsto C(-,M)$. By Theorems \ref{thm3a} and \ref{thm4a}, the category $\mathbf{CMod}(\underline{R})$ is abelian with nice properties, such as having all (small) limits and colimits and also enough projectives. So, we might ask if we can use the (co)homology in $\mathbf{CMod}(\underline{R})$ to define/understand (co)homology in $\TModR$. This depends, among other things, on whether the condensation functor $C\colon\TModR\to\mathbf{CMod}(\underline{R})$ is nice (such as whether it is fully faithful).

Recall (\cite[page 8]{condensed}) that a topological space $X$ is \emph{compactly generated} if a map $f\colon X\to Y$ to another space $Y$ is continuous if and only if for every continuous map $S\to X$ from a compact Hausdorff space $S$, the composition $S\to X\to Y$ is continuous. Since every compact Hausdorff space admits a surjection from a compact Hausdorff extremally disconnected space, and every continuous map between compact Hausdorff spaces is a quotient map, we may replace ``compact Hausdorff $S$" in the definition of compactly generated by ``profinite $S$" or ``compact Hausdorff extremally disconnected $S$" if we wish.

Note that there are some set-theoretic problems here (\cite[page 8 and Warning 2.14]{condensed}). We should either fix an uncountable strong limit cardinal $\kappa$ and consider only $\kappa$-compactly generated spaces and $\kappa$-small condensed sets, or restrict to $T1$ topological spaces (see \cite[Proposition 2.15]{condensed}). We will mostly be interested in discrete spaces and profinite spaces, which are $T1$ and compactly generated.

\begin{prop}\label{Ffull}
Let $R$ be a topological ring. The condensation functor $C\colon\TModR\to(\kappa\text{-})\mathbf{CMod}(\underline{R})\colon M\mapsto\underline{M}$ is faithful, and fully faithful when restricted to the full subcategory of all $M$ that are $(\kappa$-$)$compactly generated as topological spaces.
\begin{proof}
Exactly the same as in \cite[Proposition 1.7]{condensed}.
\end{proof}
\end{prop}

Now suppose $R$ is a profinite ring. Then the categories $\PModR$ of profinite $R$-modules and $\DModR$ of discrete (topological) $R$-modules are abelian. Often, the case of interest is when $R$ is a profinite group algebra i.e.\ $R=\kG=\varprojlim k[G_i]$, where $k$ is a commutative profinite ring and $G=\varprojlim G_i$ is a profinite group.

\begin{thm}\label{PDfullsub}
Let $R$ be a profinite ring. Then $\PModR$ and $\DModR$ are full abelian subcategories of $\mathbf{CMod}(\underline{R})$ via the condensation functor $C$.
\begin{proof}
Proposition \ref{Ffull} shows that the restriction of $C$ to either $\PModR$ or $\DModR$ is fully faithful, since profinite and discrete spaces are compactly generated. It remains to show that $C$ is exact. First note that kernels and cokernels in $\PModR$ (or $\DModR$) are the usual algebraic ones with induced topology, and that kernels and cokernels in $\mathbf{CMod}(\underline{R})$ are calculated pointwise (Theorem \ref{thm3a}).

Suppose we have an exact sequence $M\to N\overset{f}{\longrightarrow} K$ in $\PModR$ or $\DModR$. We need to show that for each $S\in\CHED$, $C(S,M)\to C(S,N)\to C(S,K)$ is exact. Suppose $g\in C(S,N)$ satisfies $f\circ g=0$. For $\DModR$, an arbitrary pointwise lift of $g$ to a function $S\to M$ is continuous. For $\PModR$, note that extremally disconnected spaces are projective in $\Pro$ (see \cite{gleason}).
\end{proof}
\end{thm}

Remark: A similar argument shows that $\PAb=\mathbf{PMod}(\widehat{\Z})$ and $\mathbf{DMod}(\widehat{\Z})$ (i.e.\ discrete torsion abelian groups) are full abelian subcategories of $\CAb$.

The categories $\PModR$ and $\DModR$ are abelian and related by Pontryagin Duality (see \cite{profinite}, for example), but either one alone feels somewhat lacking. The category $\PModR$ has all small limits which satisfy (AB4*), (AB5*) and has enough projectives, but in general the coproduct is not exact. The category $\DModR$ suffers dual problems. Now that we have embedded both $\PModR$ and $\DModR$ in a single abelian category, namely $\mathbf{CMod}(\underline{R})$, the next natural question to ask is whether their cohomologies agree in some way.

We have a few different cohomology theories at hand. Firstly, because $\PModR$, $\DModR$ and $\mathbf{CMod}(\underline{R})$ are abelian, we can define $\Ext^*(-,-)$ using Yoneda extensions. All of these categories have either enough projectives or enough injectives, so we can also take projective/injective resolutions as usual. On the other hand, we can take $\Ext^*_R(A,B)$ with $A\in\PModR$ and $B\in\DModR$ as the (classical) right derived functor of $\Hom_R(A,B)$ in either variable (i.e.\ calculated by projectively resolving $A$ or injectively resolving $B$ or both). Before we proceed, we need to talk about \emph{solid} modules.

\subsection{Solid Modules}
\label{sub2}

The following definitions are from \cite[Definition 5.1]{condensed} and \cite[Definition 2.3.25]{brink}:

\begin{defn}\label{solid}
Given a profinite space $S=\varprojlim S_i$ ($S_i$ finite), the \emph{free solid abelian group} on $S$ is the condensed abelian group $\Z[\underline{S}]^\blacksquare=\varprojlim\Z[\underline{S_i}]$. There is a natural map of condensed sets $\underline{S}\to\Z[\underline{S}]^\blacksquare$ (hence a map of condensed abelian groups $\Z[\underline{S}]\to\Z[\underline{S}]^\blacksquare$).

A condensed abelian group $\M$ is \emph{solid} if for all $S\in\Pro$ and all maps $f\colon\underline{S}\to\M$ in $\CSet$, there is a unique map $\Z[\underline{S}]^\blacksquare\to\M$ extending $f$.

A condensed ring $\R$/condensed $\R$-module $\M$ is \emph{solid} if its underlying condensed abelian group is solid.
\end{defn}

By \cite[Proposition 5.7]{condensed}, a free solid abelian group is indeed solid, justifying the name. Moreover, by \cite[Proposition 5.6]{condensed}, we could have replaced $\Pro$ in the definition above everywhere with $\CHED$ without changing the notions. In addition, \cite[Proposition 2.1]{analytic} shows that the map $\Z[\underline{S}]\to\Z[\underline{S}]^\blacksquare$ is actually injective. As explained in \cite[lectures I and II]{analytic}, the motivation for defining solid abelian groups is that they are ``complete" in some sense, and we can enlarge or ``complete" the free condensed abelian group $\Z[\underline{S}]$ to the free solid abelian group $\Z[\underline{S}]^\blacksquare$.

\begin{thm}[\cite{condensed} Theorem 5.4, Corollary 5.5]\label{freeab} For a profinite space $S$, the abelian group $C(S,\Z)$ of continuous maps is free. In particular, if $C(S,\Z)=\bigoplus_I\Z$, then $\Z[\underline{S}]^\blacksquare=\prod_I\underline{\Z}$.
\end{thm}

It is useful to know how the free solid abelian groups are isomorphic to products of $\underline{\Z}$. The argument uses internal Hom. First note that by \cite[Proposition 4.2]{condensed}, if $A$ and $B$ are Hausdorff topological abelian groups with $A$ compactly generated, then there is a natural isomorphism of condensed abelian groups $\UHom(\underline{A},\underline{B})=\underline{\Hom(A,B)}$, where $\Hom(A,B)$ has the compact-open topology. To avoid drowning in underlines, we might sometimes write $\UHom(A,B)$ for $\UHom(\underline{A},\underline{B})$ (although we have been very careful with underlines so far). Now if $C(S,\Z)=\bigoplus_I \Z$, then $$\Z[\underline{S}]^\blacksquare=\varprojlim\Z[\underline{S_i}]=\varprojlim\UHom(C(S_i,\Z),\Z)=\UHom(C(S,\Z),\Z)=\UHom\left(\bigoplus_I\Z,\Z\right)=\prod_I\underline{\Z}.$$

The isomorphism $C(S,\Z)=\bigoplus_I\Z$ is highly non-canonical, but once we fix such isomorphisms $C(S,\Z)=\bigoplus_I\Z$ and $C(T,\Z)=\bigoplus_J\Z$ for $S,T$ profinite, we canonically have $C(S\times T,\Z)=C(S,\Z)\otimes C(T,\Z)=\bigoplus_{I\times J}\Z$. (The first isomorphism can be proven by writing $S$ and $T$ as inverse limits of finite spaces.)

Let $\SAb$ denote the full subcategory of solid abelian groups in $\CAb$, and $\SR$ the full subcategory of solid $\R$-modules in $\CModRR$. The category $\SAb$ forms a subcategory of $\CAb$ with favourable properties:

\begin{thm}[\cite{condensed} Theorem 5.8, Corollary 6.1, Theorem 6.2, Proposition 6.3]\label{thmsolid}
The category $\SAb$ of solid abelian groups is an abelian subcategory  of $\CAb$ closed under all limits, colimits and extensions. The compact projective objects of $\SAb$ are exactly the solid objects $\prod_I\underline{\Z}$, for $I$ any set, and furthermore they generate. The inclusion $\SAb\inj\CAb$ has a left adjoint $\M\mapsto\M^\blacksquare$, called the \emph{solidification functor}, which is the unique colimit-preserving extension of $\Z[\underline{S}]\mapsto\Z[\underline{S}]^\blacksquare$.

There is a unique way to endow $\SAb$ with a symmetric monoidal tensor product $\botimes$ such that the solidification functor $\M\mapsto\M^\blacksquare$ is symmetric monoidal. Moreover, we have $$\left(\prod_I\underline{\Z}\right)\botimes\left(\prod_J\underline{\Z}\right)=\prod_{I\times J}\underline{\Z}.$$
\end{thm}

The solidification functor $(-)^\blacksquare$ establishes $\SAb$ as a reflective subcategory of $\CAb$, similar to sheafification being the reflector of the embedding of sheaves into presheaves. In particular, if $\M\in\CAb$ is already solid, then $\M=\M^\blacksquare$ i.e.\ the natural map $\M\to\M^\blacksquare$ is an isomorphism. Note also that for $\M\in\SAb$, $\M\botimes(-)=(\M\otimes-)^\blacksquare\colon\SAb\to\SAb$ commutes with all colimits.

Let $\R$ be a condensed ring. Recall that we initially defined an $\R$-module pointwise (Definition \ref{condmoda}), which we then showed is equivalent to an $\R$-module object in the monoidal category ($\CAb$, $\otimes$) (Proposition \ref{condmod2a}). A key property we used is that sheafification is a monoidal functor (Lemma \ref{sheaftensora}). Indeed, mimicking Proposition \ref{condmod2a}, we obtain the following proposition (see also \cite[Lemma 2.3.27]{brink}).

\begin{prop}\label{solidpr}
A solid ring $\R$ is equivalently a monoid in $\left(\SAb,\botimes\right)$. A solid $\R$-module (for $\R$ solid) is equivalently an $\mathcal{R}$-module object in $\SAb$.

If $\R\in\CRing$ and $\M\in\CModRR$, then $\R^\blacksquare$ is naturally a solid ring and $\M^\blacksquare$ is naturally an $\R^\blacksquare$-module. In particular, $\M^\blacksquare$ also has an $\R$-module structure induced by the map $\R\to\R^\blacksquare$. The solidification functor induces a functor $(-)^\blacksquare\colon\CModRR\to\SR$ (also called solidification) which is left adjoint to the inclusion $\SR\inj\CModRR$. This in turn gives a functor $$\SAb\to\SR\colon \M\mapsto\R^\blacksquare\otimes^\blacksquare\M$$ which is left adjoint to the forgetful functor $\SR\to\SAb$.

There is an equivalence of categories $\SR=\mathbf{Solid}(\R^\blacksquare)$.
\end{prop}

\begin{corr}\label{solidcor}
Let $\R$ be a condensed ring. The category $\SR$ of solid $\R$-modules is an abelian subcategory of $\CModRR$ closed under all limits, colimits and extensions. Moreover, $\SR$ is generated by the compact projective objects $(\R\otimes\Z[\underline{S}])^\blacksquare$, $S\in\CHED$. In particular, $\SR$ has enough projectives.
\begin{proof}
The first claim follows at once from Theorem \ref{thmsolid} and the fact that the forgetful functor $\CModRR\to\CAb$ commutes with taking limits and colimits (Theorem \ref{thm3a}). The other claims follow from $(-)^\blacksquare$ being left adjoint to $\SR\inj\CModRR$ and the description of compact projective generators for $\CModRR$ (Theorem \ref{thm4a}).
\end{proof}
\end{corr}

Thus, the category $\SR$ of solid $\R$-modules has basically the same favourable categorical properties as $\CModRR$.

\subsection{Solid Modules over Profinite Rings}
\label{sub3}

Given a profinite ring $R$, recall that we have the condensation functor $C\colon\TModR\to\mathbf{CMod}(\underline{R})\colon M\mapsto\underline{M}$. Let's write $C^{P}$ and $C^{D}$ for the restrictions of $C$ to $\PModR$ and $\DModR$ respectively, so $C^P$ and $C^D$ are fully faithful and exact (Theorem \ref{PDfullsub}). The reason we need solid modules will become evident below.

\begin{lemma}\label{smalllemma} Let $R$ be a profinite ring and $C^P\colon\PModR\inj\mathbf{CMod}(\underline{R})$, $C^D\colon\DModR\inj\mathbf{CMod}(\underline{R})$ be the condensation functors. Then:
\begin{enumerate}[label=(\roman*)]
\item The functor $C^P$ preserves limits.
\item The functor $C^D$ preserves colimits.
\end{enumerate}
\begin{proof}
\begin{enumerate}[label=(\roman*)]
\item Follows from the fact that limits in $\PModR$ are the algebraic ones (with induced topology) and limits in $\mathbf{CMod}(\underline{R})$ are taken pointwise (Theorem \ref{thm3a}).
\item The functor $\underline{(-)}\colon\mathbf{Ab}\to\CAb$ preserves colimits by Lemma \ref{samesheafa}. This implies the claim because $\underline{(-)}$ and arbitrary colimits all commute with finite products (which are finite coproducts in $\Ab$).
\end{enumerate}
\end{proof}
\end{lemma}

In particular, the embedding $\PAb\inj\CAb$, which can be viewed as the composition $\PAb=\mathbf{PMod}(\widehat{\Z})\to\mathbf{CMod}(\underline{\widehat{\Z}})\to\CAb$ (where the last map is the forgetful functor), also preserves limits.

\begin{lemma}[\cite{brink} Lemma 2.3.18]\label{brinklemma1} If $M$ is a limit of discrete abelian groups in $\TAb$, then $\underline{M}$ is a solid abelian group. In particular, if $R$ is a profinite ring, then $\underline{R}$ is solid, and the essential images of the condensation functors $C^P$ and $C^D$ lie in $\mathbf{Solid}(\underline{R})$.
\begin{proof}
We will include a proof for completeness. First suppose $M$ is discrete, so in particular it is the cokernel of some map $\bigoplus_i\Z\to\bigoplus_j\Z$. By Lemma \ref{samesheafa}, the functor $\underline{(-)}\colon\mathbf{Ab}\to\CAb$ is a left adjoint, so preserves all colimits. Thus, $\underline{M}$ is the cokernel of $\bigoplus_i\underline{\Z}\to\bigoplus_j\underline{\Z}$, but $\underline{\Z}$ is solid and $\SAb$ is closed under all colimits (Theorem \ref{thmsolid}), so $\underline{M}$ is solid. Now note that $\underline{(-)}\colon\TAb\to\CAb$ preserves all limits and that $\SAb$ is closed under all limits. The rest follows.
\end{proof}
\end{lemma}

Henceforth, when we write $C^P$ we will mean the embedding $\PModR\inj\mathbf{Solid}(\underline{R})$ unless otherwise specified (though sometimes it does not matter if we mean $\PModR\inj\mathbf{Solid}(\underline{R})$ or $\PModR\inj\mathbf{CMod}(\underline{R})$).

The next proposition (which should be viewed as the key proposition of the entire paper) and its importance are essentially suggested to the author by Peter Scholze, whom the author is in debt to. Any mistakes in its statement or proof are solely due to the author.

It was later pointed out to the author by Guido Bosco that the following proposition has already appeared as \cite[Lemma A.19]{guido}.

\begin{prop}\label{limtensor}
Let $M$ be a profinite abelian group. Then $\underline{M}\botimes(\prod_I\underline{\Z})=\prod_{I}\underline{M}$ as solid abelian groups. This is natural in $M$, so if in addition $R$ is a profinite ring and $M$ is an $R$-module, then the above isomorphism is automatically an $\underline{R}$-module isomorphism. In particular, $\mathbf{Solid}(\underline{R})$ is generated by the compact projective objects $\prod_I\underline{R}$ for varying sets $I$.
\begin{proof}
(cf.\ \cite[Example 6.4]{condensed}) By \cite[Proposition 4.2/Theorem 4.5]{condensed}, there is a functorial partial resolution of condensed abelian groups $\Z[\underline{M^2}]\to\Z[\underline{M}]\to\underline{M}\to0$. Solidifying this gives a partial resolution of solid abelian groups $\Z[\underline{M^2}]^\blacksquare\to\Z[\underline{M}]^\blacksquare\to\underline{M}\to0$. By Theorems \ref{freeab} and \ref{thmsolid}, if $C(M,\Z)=\bigoplus_J\Z$, then the resolution can be written as $\prod_{J\times J}\underline{\Z}\to\prod_{J}\underline{\Z}\to\underline{M}\to0$. Using the exactness of products (AB4*) and the last part of Theorem \ref{thmsolid} again, we see that both $\underline{M}\botimes(\prod_I\underline{\Z})$ and $\prod_{I}\underline{M}$ are cokernels of the induced map $\prod_{I\times J\times J}\underline{\Z}\to\prod_{I\times J}\underline{\Z}$. Naturality in $M$ follows from the functoriality of the original partial resolution.

To show that the isomorphism is automatically an $\underline{R}$-module isomorphism if $M$ is an $R$-module, we need to talk about the profinite tensor product $\cotimes$, which we postpone; see Remark \ref{rmkpf} for the proof. Alternatively, we can prove this directly as follows. Note, by mimicking the proof above, that $\underline{R}\botimes(\prod_I\underline{M})=\prod_I(\underline{R}\botimes\underline{M})$ as solid abelian groups. We can then show that the diagram
\[\begin{tikzcd}
\underline{R}\botimes\underline{M}\botimes(\prod_I\underline{\Z}) \arrow[d] \arrow[r] & \underline{R}\botimes(\prod_{I}\underline{M}) \arrow[d] \\
\underline{M}\botimes(\prod_I\underline{\Z}) \arrow[r]                                & \prod_{I}\underline{M}                               
\end{tikzcd}\]
commutes by checking that it commutes on the level of presheaves, which is exactly what we need. (Note that here the bottom map is induced by the presheaf map $m\otimes(n_i)\mapsto(n_im)$.)

The last part follows from Corollary \ref{solidcor}.
\end{proof}
\end{prop}

As we will be using the above proposition quite often, let us make more transparent what we have just done. Consider the special case of the above when $\prod_I\underline{\Z}=\Z[\underline{S}]^\blacksquare$ is free with $C(S,\Z)=\bigoplus_I\Z$. We can apply the functor $(-)\botimes\Z[\underline{S}]^\blacksquare$ to the partial resolution of $\underline{M}$ to get $$\Z[\underline{M^2}]^\blacksquare\botimes\Z[\underline{S}]^\blacksquare\to
\Z[\underline{M}]^\blacksquare\botimes\Z[\underline{S}]^\blacksquare\to\underline{M}\botimes\Z[\underline{S}]^\blacksquare\to0,$$ which can also be written as \begin{eqnarray*}\UHom(C(M^2,\Z),\Z)\botimes\UHom(C(S,\Z),\Z)&\to&\UHom(C(M,\Z),\Z)\botimes\UHom(C(S,\Z),\Z)\\&\to&\underline{M}\botimes\UHom(C(S,\Z),\Z)\to0.\end{eqnarray*} On the other hand, we can also apply the (exact!) functor $\prod_I(-)=\UHom(\underline{C(S,\Z)},-)$ to get $$\UHom(\underline{C(S,\Z)},\Z[\underline{M^2}]^\blacksquare)\to\UHom(\underline{C(S,\Z)},\Z[\underline{M}]^\blacksquare)\to\UHom(\underline{C(S,\Z)},\underline{M})\to0,$$ which, by the Hom-tensor adjunction (or rather its internal version), can be written as $$\UHom(C(M^2\times S,\Z),\Z)\to\UHom(C(M\times S,\Z),\Z)\to\UHom(C(S,\Z),M)\to0.$$ The first two terms of these two partial resolutions can be identified exactly as in the proof of \cite[Proposition 6.3]{condensed}, which proves the above proposition in this special case. The general case then follows by writing $\prod_I\underline{\Z}$ as a retract of a free solid abelian group.

For $\R$ a condensed ring, there is an obvious alternative notion of solid $\R$-modules analogous to the first part of Definition \ref{solid}. Namely, given a profinite space $S=\varprojlim S_i$, we define the condensed $\R$-module $\R[\underline{S}]^\square=\varprojlim(\R^\blacksquare\otimes\Z[\underline{S_i}])$ (we are not using $\R[\underline{S}]^\blacksquare$ because that will later mean $(\R\otimes\Z[\underline{S}])^\blacksquare$ when we define condensed group rings). Note that this is solid as a condensed abelian group. There is a natural map of condensed sets $\underline{S}\to\R[\underline{S}]^\square$ (hence a map of $\R$-modules $\R\otimes\Z[\underline{S}]\to\R[\underline{S}]^\square$). We can then call an $\R$-module $\M$ \emph{alternate solid} if for all $S\in\Pro$ and all maps $f\colon\underline{S}\to\M$ in $\CSet$, there is a unique $\R$-module map $\R[\underline{S}]^\square\to\M$ extending $f$. Note that this is the definition of solid $\underline{\F_p}$-modules given on \cite[page 14]{analytic}, where it is also noted that this is equivalent to being solid as a condensed abelian group. Indeed, as a consequence of Proposition \ref{limtensor}:

\begin{corr}\label{samesoliddef}
Let $R$ be a profinite ring. Then alternate solid $\underline{R}$-modules are precisely solid $\underline{R}$-modules.
\begin{proof}
An $\UR$-module $\M$ is solid if it is solid as a condensed abelian group i.e.\ if the natural map $\Hom(\Z[\underline{S}]^\blacksquare,\M)\to\Hom(\Z[\underline{S}],\M)$ is an isomorphism for each $S\in\Pro$. On the other hand, $\M$ is alternate solid if the natural map $\Hom_{\underline{R}}\left(\varprojlim(\underline{R}\otimes\Z[\underline{S_i}]),\M\right)\to\Hom_{\underline{R}}(\underline{R}\otimes\Z[\underline{S}],\M)$ is an isomorphism for each $S\in\Pro$. Now note that by Proposition \ref{limtensor}, if $C(S,\Z)=\bigoplus_I\Z$, then \begin{eqnarray}\label{eqntensor}\varprojlim(\underline{R}\otimes\Z[\underline{S_i}])=\varprojlim\UHom(C(S_i,\Z),R)=\UHom(C(S,\Z),R)=\prod_I\underline{R}=(\underline{R}\otimes\Z[\underline{S}])^\blacksquare.\end{eqnarray} It is thus clear that if $\M$ is solid then it is alternate solid.

The other direction is slightly trickier. The following argument is essentially the same as the first part of \cite[Lemma 5.9]{condensed}, but we include it for the convenience of the reader and without referring to derived categories. Suppose $\M$ is alternate solid. We can find a surjection of $\UR$-modules $\bigoplus(\underline{R}\otimes\Z[\underline{S_j}])\sur\M$ which, by alternate solidity, gives a surjection $\bigoplus(\underline{R}\otimes\Z[\underline{S_j}])^\blacksquare\sur\M$. Write $\mathcal{L}=\bigoplus(\underline{R}\otimes\Z[\underline{S_j}])^\blacksquare$ and let $\K$ be the kernel of $\mathcal{L}\sur\M$. Note that being alternate solid is closed under kernels by the 5-lemma, and that $\mathcal{L}$ is alternate solid because it is solid. Hence, $\K$ is alternate solid, and by repeating the above argument, we find a partial resolution $\mathcal{L}'\to\mathcal{L}\sur\M$ with $\mathcal{L}'$, $\mathcal{L}$ solid, so $\M$ is also solid.
\end{proof}
\end{corr}

\subsection{Preservation of $\Hom$ and $\Ext$}
\label{sub4}

We are almost ready to prove that solid cohomology extends profinite cohomology. Let $R$ be a profinite ring. Since $\PModR$ is a full subcategory of $\mathbf{CMod}(\underline{R})$ (or $\mathbf{Solid}(\underline{R})$), we clearly have $\Hom_{\underline{R}}(\underline{M},\underline{N})=\Hom_R(M,N)$ for $M,N$ profinite $R$-modules, but this is even true on the condensed level:

\begin{prop}\label{homunderline}
Let $R$ be a topological ring, and $M, N$ be Hausdorff topological $R$-modules with $M$ compactly generated. Then there is a natural isomorphism of condensed abelian groups $$\UHom_{\underline{R}}(\underline{M},\underline{N})=\underline{\Hom_R(M,N)},$$ where $\Hom_R(M,N)$ has the compact-open topology.
\begin{proof}
This is similar to \cite[Proposition 4.2]{condensed}, which is the above statement for $R=\Z$. The quoted proposition constructs an explicit (well-defined) map \begin{eqnarray}\label{eqn1}\Hom(\underline{M}\otimes\Z[\underline{S}],\underline{N})=\UHom(\underline{M},\underline{N})(S)\to\underline{\Hom(M,N)}(S)=C(S,\Hom(M,N))\end{eqnarray} which evaluates the underlying abelian groups of the left hand side, and then shows that the map is an isomorphism. Note that if we take a map in $\UHom_{\underline{R}}(\underline{M},\underline{N})(S)=\Hom_{\underline{R}}(\underline{M}\otimes\Z[\underline{S}],\underline{N})$, then its image under (\ref{eqn1}) will send any $s\in S$ to an $R$-homomorphism $M\to N$. This shows that (\ref{eqn1}) induces a well-defined injection $\UHom_{\underline{R}}(\underline{M},\underline{N})(S)\inj\underline{\Hom_R(M,N)}(S)$. It only remains to prove that this map is also surjective, which one can do without difficulty by following the proof of \cite[Proposition 4.2]{condensed}.
\end{proof}
\end{prop}

\begin{lemma}\label{solidhoml}
Let $\R$ be a condensed ring and $\M,\mathcal{N}\in\CModRR$ with $\mathcal{N}$ solid. Then their enriched Hom $\UHom_\R(\M,\mathcal{N})\in\CAb$ is also solid.
\begin{proof}
Follows from the Hom-tensor adjunction (Proposition \ref{homtensora}):
\begin{eqnarray*}
\Hom(\Z[\underline{S}]^\blacksquare,\UHom_{\R}(\M,\mathcal{N}))&=&\Hom_{\R}(\Z[\underline{S}]^\blacksquare\otimes\M,\mathcal{N})\\
&=&\Hom_{\R}(\Z[\underline{S}]^\blacksquare\botimes\M^\blacksquare,\mathcal{N}) \\
&=&\Hom_{\R}(\Z[\underline{S}]\otimes\M,\mathcal{N})\\
&=&\Hom(\Z[\underline{S}],\UHom_{\R}(\M,\mathcal{N})).
\end{eqnarray*}
\end{proof}
\end{lemma}

For a solid $\R$-module $\mathcal{N}$, we write $\UExt^*_{\SR}(-,\mathcal{N})$ for the (classical) right derived functor of $\UHom_{\R}(-,\mathcal{N})\colon\SR^{\mathrm{op}}\to\SAb$. Note that $\UExt^*_{\SR}(-,\mathcal{N})(*)=\Ext^*_{\SR}(-,\mathcal{N})$ is the right derived functor of $\Hom_{\R}(-,\mathcal{N})\colon\SR^{\mathrm{op}}\to\Ab$.

Let $R$ be a profinite ring and $M,N\in\PModR$ with $N=\varprojlim N_i$ ($N_i$ finite). Note that if $M$ is (topologically) finitely generated over $R$, then $\Hom_R(M,N)=\varprojlim\Hom_R(M,N_i)$ is profinite when endowed with the compact-open topology. Recall that $M$ is \emph{of type $FP_\infty$} if it has a projective resolution $P$ with each term finitely generated. In this case, each $\Ext^i_R(M,N)=H^i(\Hom_R(P,N))$ is naturally a profinite abelian group. See \cite[Lemma 2.2.17 and Section 3.1]{cookthesis} for more details.

\begin{thm}\label{thm5}
Let $R\in\PRing$ and let $C^P\colon\PModR\inj\mathbf{Solid}(\underline{R})$ be the (fully faithful and exact) condensation functor. Then:
\begin{enumerate}[label=(\roman*)]
\item The functor $C^P$ preserves projectives.
\item Let $M,N$ be profinite $R$-modules. Then $$\Ext^*_{\mathbf{Solid}(\underline{R})}(\underline{M},\underline{N})=\Ext^*_R(M,N).$$ If moreover $M$ is of type $FP_\infty$, then $$\UExt^*_{\mathbf{Solid}(\underline{R})}(\underline{M},\underline{N})=\underline{\Ext^*_R(M,N)},$$ where the $\Ext^*_R(M,N)$ have their natural profinite topologies.
\item Let $M$ be a profinite $R$-module and $A$ be a discrete $R$-module. Then $$\UExt^*_{\mathbf{Solid}(\underline{R})}(\underline{M},\underline{A})=\underline{\Ext^*_R(M,A)},$$ where the $\Ext^*_R(M,A)$ have discrete topologies.
\end{enumerate}
\begin{proof}
\begin{enumerate}[label=(\roman*)]
\item As every projective profinite $R$-module is a direct summand of a free one, it suffices to show that $C^P$ sends free $R$-modules to projective solid $\underline{R}$-modules. Let $S=\varprojlim S_i$ be a profinite space and $R\llbracket S\rrbracket=\varprojlim R[S_i]$ be the corresponding free profinite module. Then equation (\ref{eqntensor}) from the proof of Corollary \ref{samesoliddef} shows that $C^P(R\llbracket S\rrbracket)=(\underline{R}\otimes\Z[\underline{S}])^\blacksquare$, so we are done by Proposition \ref{limtensor}.
\item The first part follows from (i) by taking projective resolutions, while the second part follows from (i), Proposition \ref{homunderline} and the paragraph above the theorem.
\item Follows from (i) and Proposition \ref{homunderline} by taking projective resolutions. Note that for $M$ profinite and $A$ discrete, the compact-open topology on $\Hom_R(M,A)$ is discrete.
\end{enumerate}
\end{proof}
\end{thm}

\begin{exmp}
We can of course take $R=\widehat{\Z}$ in the theorem above, so that the embedding $\PAb\inj\mathbf{Solid}(\underline{\widehat{\Z}})$ preserves projectives. However, it is not true that the embedding $\PAb\inj\SAb$ preserves projectives. Indeed, $\widehat{\Z}$ is a projective profinite abelian group, but $\underline{\widehat{\Z}}$ is not projective in $\SAb$. To see this, note that there is a surjection $\prod_I\underline{\Z}\sur\underline{\widehat{\Z}}$ of solid abelian groups, and that $\Hom(\underline{\widehat{\Z}},\underline{\Z})=\Hom(\widehat{\Z},\Z)=0$.
\end{exmp}

Thus, solid cohomology extends profinite cohomology in the sense above. In fact, as we will see below, the functor $C^P$ also preserves tensor products.

\subsection{Preservation of Tensor Products}
\label{sub5}

\begin{prop}\label{invtensor}
Let $M=\varprojlim M_i$ and $N=\varprojlim N_j$ be inverse limits of finite abelian groups $M_i$ and $N_j$ respectively. Then $\underline{M}\botimes\underline{N}=\varprojlim_{i,j}\left(\underline{M_i}\otimes \underline{N_j}\right)$.
\begin{proof}
As in Proposition \ref{limtensor}, we take, for each $i$, the partial resolution of abelian groups \begin{eqnarray}\label{eqn2}\Z[M_i^2]\to\Z[M_i]\to M_i\to0.\end{eqnarray} Applying the condensation functor and taking inverse limits recovers the resolution $\Z[\underline{M^2}]^\blacksquare\to\Z[\underline{M}]^\blacksquare\to\underline{M}\to0$ of $\underline{M}$, so $\underline{M}\botimes\underline{N}$ is the cokernel of the map $(\Z[\underline{M^2}]\otimes\underline{N})^\blacksquare\to(\Z[\underline{M}]\otimes\underline{N})^\blacksquare$, where (on the level of presheaves) a generator $[(a,b)]$ of $\Z[\underline{M^2}]$ is sent to $[a+b]-[a]-[b]$ of $\Z[\underline{M}]$. On the other hand, we can either apply $(-)\otimes N_j$ to (\ref{eqn2}) and then condense it, or condense (\ref{eqn2}) and then apply $(-)\otimes\underline{N_j}$. This shows that $\underline{M_i}\otimes \underline{N_j}=\underline{M_i\otimes N_j}$, so $\underline{M_i}\otimes \underline{N_j}$ is already solid. Finally, applying $(-)\otimes N_j$ to (\ref{eqn2}), taking inverse limits and condensing shows that $\varprojlim\underline{M_i}\otimes \underline{N_j}$ is the cokernel of the map $\UHom(C(M^2,\Z),N)\to\UHom(C(M,\Z),N)$, which can be identified with the map $(\Z[\underline{M^2}]\otimes\underline{N})^\blacksquare\to(\Z[\underline{M}]\otimes\underline{N})^\blacksquare$ through Proposition \ref{limtensor}. We used here that inverse limits of profinite abelian groups are exact (even though inverse limits in $\CAb$ might not be in general).
\end{proof}
\end{prop}

Recall that there is a completed tensor product $\cotimes=\cotimes_{\widehat{\Z}}\,$ of profinite abelian groups defined, for $M=\varprojlim M_i$ and $N=\varprojlim N_j$ ($M_i$, $N_j$ finite), by $M\cotimes N=\varprojlim M_i\otimes_{\widehat{\Z}} N_j=\varprojlim M_i\otimes_{\Z} N_j$. As an immediate consequence of the previous proposition, we have:

\begin{corr}\label{corrtensor}
Let $M$ and $N$ be profinite abelian groups. Then $\underline{M}\botimes\underline{N}=\underline{M\cotimes N}$ as solid abelian groups. In particular, if $M=\varprojlim M_i$ and $N=\varprojlim N_j$ are inverse limits of profinite abelian groups $M_i$ and $N_j$ (not necessarily finite), then $\underline{M}\botimes\underline{N}=\varprojlim_{i,j}\left(\underline{M_i}\botimes \underline{N_j}\right)$. As a consequence, $(\prod_I\underline{M})\botimes(\prod_J\underline{N})=\prod_{I\times J}(\underline{M}\botimes\underline{N})$.
\begin{proof}
The first part is immediate from above, while for the next part, note that the analogous statement $M\cotimes N=\varprojlim M_i\cotimes N_j$ for profinite abelian groups is true. (We could have also proven the final claim earlier and more directly by mimicking the proof of Proposition \ref{limtensor}.)
\end{proof}
\end{corr}

\begin{rmk}\label{rmkpf}
In Proposition \ref{limtensor} (which states that $\underline{M}\botimes(\prod_I\underline{\Z})=\prod_{I}\underline{M}$ if $M$ is a profinite abelian group, with the isomorphism being natural in $M$), we claimed that if $R$ is a profinite ring and $M$ is also an $R$-module, then the above map is automatically an $\underline{R}$-module map. Let us prove this now. First note that we do not need this additional fact to prove the above corollary.

The $R$-module structure on $M$ can be described by a profinite abelian group map $R\cotimes M\to M$. Applying the naturality of the isomorphism  $\underline{M}\botimes(\prod_I\underline{\Z})=\prod_{I}\underline{M}$ to this map and using the above corollary, we see that the diagram
\[\begin{tikzcd}
\underline{R}\botimes\underline{M}\botimes(\prod_I\underline{\Z}) \arrow[d] \arrow[r] & \prod_{I}(\underline{R}\botimes\underline{M}) \arrow[d] \\
\underline{M}\botimes(\prod_I\underline{\Z}) \arrow[r]                                & \prod_{I}\underline{M}                                 
\end{tikzcd}\]
commutes, which is what we want.
\end{rmk}

\begin{exmp}
Let $p$ and $q$ be distinct primes. Then $$\underline{\Z_p}\botimes\underline{\Z_p}=\underline{\varprojlim(\Z/p^i\Z\otimes\Z/p^j\Z)}=\underline{\Z_p}$$ and $$\underline{\Z_p}\botimes\underline{\Z_q}=\underline{\varprojlim(\Z/p^i\Z\otimes\Z/q^j\Z)}=0.$$ (cf.\ \cite[Example 6.4]{condensed})
\end{exmp}

Recall from Section \ref{sec1} that given a condensed ring $\R$ and two $\R$-modules $\M$ and $\mathcal{N}$ (which are appropriately sided), we have their tensor product $\M\otimes_\R\mathcal{N}\in\CAb$.\label{solidt} This descends to a solid tensor product $\brotimes$ by defining, for $\M,\mathcal{N}\in\SR$, $\M\brotimes\mathcal{N}=(\M\otimes_\R\mathcal{N})^\blacksquare\in\SAb$. Note that the solidification functor is ``monoidal" here, i.e.\ if $\M,\mathcal{N}\in\CModRR$ then $(\M\otimes_{\R}\mathcal{N})^\blacksquare=(\M^\blacksquare\otimes_{\R}\mathcal{N}^\blacksquare)^\blacksquare$. To see this, it suffices to prove the isomorphism $(\M\otimes_{\R}\mathcal{N})^\blacksquare=(\M^\blacksquare\otimes_{\R}\mathcal{N})^\blacksquare$ for free objects in $\CModRR$ and then pass to cokernels i.e.\ we want to show that $(\Z[\underline{S}]\otimes\R\otimes\Z[\underline{T}])^\blacksquare=((\Z[\underline{S}]\otimes\R)^\blacksquare\otimes\Z[\underline{T}])^\blacksquare$, but this is true because the solidification functor on $\CAb$ is monoidal (Theorem \ref{thmsolid}). One can now prove the following statements analogous to the ones above:

\begin{prop}\label{rtensors}
Let $R$ be a profinite ring and $M=\varprojlim M_i$, $N=\varprojlim N_j$ be (appropriately sided) profinite $R$-modules. Then we have the following canonical isomorphisms:
\begin{enumerate}[label=(\roman*)]
\item $\underline{M}\otimes^\blacksquare_{\underline{R}}(\prod_I\underline{R})=\prod_{I}\underline{M}$ (as solid right $\underline{R}$-modules).
\item $\underline{M}\otimes^\blacksquare_{\underline{R}}\underline{N}=\varprojlim_{i,j}\left(\underline{M_i}\otimes^\blacksquare_{\underline{R}}\underline{N_j}\right)$ (as solid abelian groups, where we can remove the $``\blacksquare"$ on the right if the $M_i$ and $N_j$ are finite).
\item $\underline{M}\otimes^\blacksquare_{\underline{R}}\underline{N}=\underline{M\crotimes N}$ (as solid abelian groups).
\item $(\prod_I\underline{M})\otimes^\blacksquare_{\underline{R}}(\prod_J\underline{N})=\prod_{I\times J}(\underline{M}\otimes^\blacksquare_{\underline{R}}\underline{N})$ (as solid abelian groups).
\end{enumerate}
\begin{proof}
Here $\crotimes$ is the completed tensor product over $R$. Part (i) follows directly from Proposition \ref{limtensor}. It only remains to prove (ii) in the case when the $M_i$ and $N_j$ are finite, along with the fact that $\underline{M_i}\otimes_{\underline{R}}\underline{N_j}=\underline{M_i\otimes_RN_j}$ in this case, since that proves (iii), which in turn proves (ii) in full. Part (iv) then follows as a consequence of (ii). To this end, let us view $M_i\otimes_RN_j$ as the cokernel of the (profinite) abelian group map \begin{eqnarray}\label{eqn3}M_i\otimes R\otimes N_j\to M_i\otimes N_j,\hspace{0.15cm} m\otimes r\otimes n\mapsto mr\otimes n-m\otimes rn.\end{eqnarray} Note that the proof of Proposition \ref{invtensor} actually shows that if one of $M$ and $N$ is finite and the other is profinite, then $\underline{M}\otimes\underline{N}=\underline{M\cotimes N}=\underline{M\otimes N}$. Thus, condensing (\ref{eqn3}) gives the natural partial resolution $\underline{M_i}\otimes \underline{R}\otimes \underline{N_j}\to \underline{M_i}\otimes \underline{N_j}\to\underline{M_i\otimes_RN_j}\to0$. On the other hand, we have by definition $\left(\underline{M_i}\otimes_{\underline{R}}\underline{N_j}\right)(S)=\underline{M_i}(S)\otimes_{\underline{R}(S)}\underline{N_j}(S)$, which is the cokernel of an abelian group map similar to (\ref{eqn3}). This gives a resolution of presheaves $\underline{M_i}\otimes \underline{R}\otimes \underline{N_j}\to \underline{M_i}\otimes \underline{N_j}\to\underline{M_i}\otimes_{\underline{R}}\underline{N_j}\to0$ which sheafifies to show that $\underline{M_i}\otimes_{\underline{R}}\underline{N_j}=\underline{M_i\otimes_RN_j}$. Finally, we can reduce the claim of (ii) (when the factors are finite) to Proposition \ref{invtensor}.
\end{proof}
\end{prop}

Remark: the reader might find it strange that we have stated the above results only for tensor products and not for $\Tor$ (like our results for $\Ext$). We can of course define and obtain results for $\Tor$ analogous to Theorem \ref{thm5}, but do note that $\Tor_*^{\mathbf{Solid}(\R)}(\M,\Nc)$ in general has to be computed by projectively resolving both $\M$ and $\Nc$, rather than just a single variable. Unlike the abstract or profinite case, projective solid abelian groups need not be flat with respect to $\botimes$. A counterexample due to Efimov (yet to appear) shows that $\prod_I{\underline{\Z}}$ is not flat when $|I|=2^{2^{\aleph_0}}$.

\subsection{Condensed Group Rings}
\label{sub6}

So far, we have only talked about condensed modules over a condensed ring. There is an obvious way to define condensed modules over a condensed group. Let $\R$ be a condensed ring and $\G$ a condensed group. If $\M$ is both an $\R$-module and a $\G$-module, we say that the actions \emph{commute} if they commute for each $S\in\CHED$. In this case, for each $S$ we can view the abelian group $\M(S)$ as a module over the abstract group ring $\R(S)[\G(S)]$, so $\M$ becomes a module over the \emph{condensed group ring} $\R[\G]=\R\otimes\Z[\G]$ (which can be viewed as the sheafification of $S\mapsto\R(S)[\G(S)]$). Everything above clearly can be reversed, so we see that $\R[\G]$-modules are exactly $\R$-modules with a commuting $\G$-action (see also \cite[page 3, point (2)]{solidcoh}).

Now let us specialise to the case where $\R=\underline{R}$ and $\G=\underline{G}$ for a commutative profinite ring $R$ and a profinite group $G$. Recall that we have the completed group algebra defined by $R\llbracket G\rrbracket=\varprojlim R[G_i]$, where $G=\varprojlim G_i$ is written as an inverse limit of finite groups. We can canonically identify $\underline{R\llbracket G\rrbracket}$ with $\underline{R}[\underline{G}]^\blacksquare$ as solid rings (see also \cite[Lemma B.4]{guido}). Indeed, it is clear that they are isomorphic as $\underline{R}$-modules (see the proof of Theorem \ref{thm5}(i)), and then we can show that their multiplicative structures agree since both are induced from the multiplications of $R$ and $G$. To be precise, we mean that the diagram
\[\begin{tikzcd}
{\underline{R}[\underline{G}]^\blacksquare\botimes\underline{R}[\underline{G}]^\blacksquare} \arrow[d] \arrow[r] & \underline{R\llbracket G\rrbracket}\botimes\underline{R\llbracket G\rrbracket} \arrow[d] \\
{\underline{R}[\underline{G}]^\blacksquare} \arrow[r]                                                            & \underline{R\llbracket G\rrbracket}                                                     
\end{tikzcd}\]
commutes, or equivalently that the diagram
\[\begin{tikzcd}
{\underline{R}\botimes\underline{R}\botimes\Z[\underline{G\times G}]^\blacksquare} \arrow[d] \arrow[r] & \underline{R\cotimes R\cotimes \widehat{\Z}\llbracket G\times G\rrbracket} \arrow[d] \\
{\underline{R}\botimes\Z[\underline{G}]^\blacksquare} \arrow[r]                                        & \underline{R\cotimes \widehat{\Z}\llbracket G\rrbracket}                            
\end{tikzcd}\]
commutes (cf.\ Remark \ref{rmkpf}). Combining everything here with the last statement of Proposition \ref{solidpr}, we obtain:

\begin{prop}\label{gmodsame}
Let $R$ be a commutative profinite ring and $G$ be a profinite group. Then the category of solid $\underline{R}$-modules with a commuting $\underline{G}$-action can be canonically identified with the category $\mathbf{Solid}(\underline{R\llbracket G\rrbracket})$ of solid $\underline{R\llbracket G\rrbracket}$-modules.
\end{prop}

Hence, if we are only interested in solid modules, then the general study of group cohomology over $\underline{G}$ in $\underline{R}$-modules is already covered by what we have been doing with ring cohomology and does not require separate treatment.

\begin{exmp}
Let $R$ be a commutative profinite ring (which we think of as being fixed) and $G$ a profinite group. For a discrete $R\llbracket G\rrbracket$-module $A$, let us write $H^*(G,A)=\Ext^*_{R\llbracket G\rrbracket}(R,A)$ for the standard cohomology groups, where $G$ acts on $R$ trivially. Note that this can be calculated using continuous group cohomology (see \cite[Theorem 6.2.4]{profinite}). We also write $\underline{H}^*_{\mathbf{Solid}}(\underline{G},\underline{A})=\UExt^*_{\mathbf{Solid}(\underline{R}[\underline{G}])}(\underline{R},\underline{A})$. Then Theorem \ref{thm5}(iii) immediately tells us that $\underline{H}^*_{\mathbf{Solid}}(\underline{G},\underline{A})=\underline{H^*(G,A)}$ (cf.\ \cite[Lemmas 2.1 and 2.5]{solidcoh}).
\end{exmp}

At this point, it seems like the concept of solidity is a (vast) generalisation of being profinite, in the sense that solid cohomology generalises certain profinite cohomology (Theorem \ref{thm5}), that solid tensor products generalise profinite tensor products (Corollary \ref{corrtensor}, Proposition \ref{rtensors}(iii)), and that solid group rings generalise profinite group rings (Proposition \ref{gmodsame}). However, we have to be careful. We should think of solidification as a kind of generalisation of profinite completion, but they are not the same. Indeed, $\underline{\Z}$ and $\underline{\widehat{\Z}}$ are very different as solid abelian groups (the underlying abelian groups are different).

\subsection{Recognising Profinite Modules}
\label{sub7}

Recall that one of our initial motivations for studying solid modules is that the category $\PModR$ of profinite modules over a profinite ring $R$ is lacking, say it does not have exact coproducts in general. We have now embedded it into the category $\mathbf{Solid}(\underline{R})$ which does have exact coproducts, and we have also shown that the classical homological notions (Ext and tensor products) extend. One natural and important question to ask is the following: can we recognise when a condensed/solid module is profinite? For example, we might start with some profinite modules and take their coproduct in $\mathbf{Solid}(\underline{R})$, and it's useful to know whether the resulting solid module lands back in $\PModR$. First observe that $\PModR$ is closed under taking limits, cokernels and extensions in $\mathbf{Solid}(\underline{R})$ (Lemma \ref{smalllemma}, Theorem \ref{PDfullsub}, Theorem \ref{thm5}(ii) respectively), so infinite coproducts are in some sense the only problem. Unfortunately, there is an issue:

\begin{exmp}
Let $R$ be a profinite ring and $M_i$, $i\in I$ be non-zero profinite $R$-modules, where $I$ is an infinite index set. Then $\bigoplus_i\underline{M_i}$, the coproduct taken in $\mathbf{Solid}(\underline{R})$, is never profinite (i.e.\ is not in the essential image of $C^P$). If it is profinite, then viewing it as a functor $\Pro^{\mathrm{op}}\to\Ab$, it must turn arbitrary colimits in $\Pro$ (which is cocomplete, see \cite[Corollary 1.4]{cookpro}) into limits in $\Ab$. Let $S$ be the coproduct of countably infinitely many points in $\Pro$ and note that $\bigoplus_i\underline{M_i}(S)=\bigoplus_i\prod_{\N}M_i$. The canonical map $\bigoplus_i\prod_{\N}M_i\to\prod_{\N}\bigoplus_iM_i$ is not surjective if $I$ is infinite and the $M_i$ are not zero, a contradiction. (Alternatively, we can simply note that the underlying abelian group of the coproduct $\coprod_iM_i$ in $\PModR$ is not the direct sum $\bigoplus_iM_i$.)
\end{exmp}

We can actually check if a condensed/solid $\UR$-module is profinite by just checking if its underlying presheaf of sets is representable, as the following shows.

\begin{lemma}
Let $R$ be a profinite ring and $\M$ be a condensed $\underline{R}$-module. Then $\M$ is profinite if and only if its underlying presheaf of sets is representable.
\begin{proof}
We have a commutative diagram
\[\begin{tikzcd}
\PModR \arrow[d] \arrow[r, "C^P", hook] & \mathbf{CMod}(\UR) \arrow[d] \\
\Pro \arrow[r, hook]                    & {[\Pro^\mathrm{op},\Set]}     
\end{tikzcd}\]
The ``only if" direction is obvious. Conversely, suppose the underlying condensed set of $\M$ is represented by $S\in\Pro$. There are maps $\M\times\M\to\M$ and $\underline{R}\times\M\to\M$ of condensed sets (defining the abelian group and module structures on $\M$), which by fully faithfulness of the Yoneda embedding give maps $S\times S\to S$ and $R\times S\to S$ of profinite spaces. It is clear that this makes $S$ a profinite $R$-module such that $\M=\underline{S}$ as condensed modules.
\end{proof}
\end{lemma}

\section{Analytic Rings}
\label{sec3}

The main goal in this short section is to show that profinite rings are analytic, which is actually an easy corollary of what we have already done.

Since $\CAb$ is an abelian category, we can consider its derived category $D(\CAb)$, which is compactly generated as a triangulated category (see \cite[page 13]{condensed}). As $\CAb$ has enough projectives and exact coproducts, $K$-projective resolutions exist (also called homotopy projective resolutions, see \cite[Proposition 4.3.4]{krause}, for example). We can thus define $-\otimes^L-$ to be the left derived functor of $-\otimes-$ and $\RUHom(-,-)$ to be the right derived functor of $\UHom(-,-)$, which can be computed using $K$-projective resolutions. These satisfy the adjunction $\Hom(\mathcal{K},\RUHom(\M,\mathcal{N}))=\Hom(\mathcal{K}\otimes^L\M,\mathcal{N})$. To prove this, note that both sides can be calculated by projectively resolving $\K$ and $\M$.

For $\R$ a condensed ring, we can similarly define its derived category of modules $D(\CModRR)$ and the derived functors $-\otimes^L_{\R}-$ and $\RUHom_{\R}(-,-)$. Recall Theorem \ref{thmsolid}, which states that the category $\SAb$ of solid abelian groups forms a nice subcategory of $\CAb$. This justifies calling the condensed ring $\underline{\Z}$ (technically $\Z_{\blacksquare}$) \emph{analytic}:

\begin{defn}[\cite{condensed} Definitions 7.1, 7.4]
A \emph{pre-analytic ring} $\A$ is a condensed ring $\underline{\A}$ together with a functor $$\CHED\to\mathbf{CMod}(\underline{\A})\colon S\mapsto\A[\underline{S}]$$ sending finite disjoint unions to products, and a natural transformation $\underline{S}\to\A[\underline{S}]$ for each $S$.

A pre-analytic ring $\A$ is \emph{analytic} if for any complex $\C\colon\ldots\to \C_1\to \C_0\to0$ of $\underline{\A}$-modules, where each $\C_i$ is a direct sum of objects of the form $\A[\underline{T}]$, the map $$\RUHom_{\underline{\A}}(\A[\underline{S}],\C)\to\RUHom_{\underline{\A}}(\underline{\A}[\underline{S}],\C)$$ in $D(\CAb)$ is an isomorphism for every $S\in\CHED$.
\end{defn}

Note that here $\underline{\A}[\underline{S}]$ means $\underline{\A}\otimes\Z[\underline{S}]$. There are a few examples of (pre-)analytic rings given in \cite[Examples 7.3]{condensed} and we would like to highlight some of them.

\begin{exmp}
\begin{enumerate}[label=(\roman*)]
\item The pre-analytic ring $\Z_\blacksquare$ has underlying condensed ring $\underline{\Z_\blacksquare}=\underline{\Z}$ with the functor $$\CHED\to\CAb\colon S\mapsto\Z_\blacksquare[\underline{S}]=\Z[\underline{S}]^\blacksquare.$$ The natural transformations $\underline{S}\to\Z[\underline{S}]^\blacksquare$ are precisely the ones used to define solid abelian groups.
\item Let $A$ be a discrete ring. Then the pre-analytic ring $(A,\Z)_\blacksquare$ has underlying condensed ring $\underline{A}$ with the functor $$\CHED\to\mathbf{CMod}(\underline{A})\colon S\mapsto\underline{A}\otimes\Z_\blacksquare[\underline{S}].$$
\item Let $R$ be a profinite ring. Then the pre-analytic ring $R_\blacksquare$ has underlying condensed ring $\underline{R}$ with the functor $$\CHED\to\mathbf{CMod}(\underline{R})\colon S\mapsto \underline{R}[\underline{S}]^\square=\varprojlim \underline{R}[\underline{S_i}].$$ Note that this is precisely what we used to define alternate solid modules (see Corollary \ref{samesoliddef}). The quoted corollary shows in particular that $\underline{R}[\underline{S}]^\square=\underline{R}[\underline{S}]^\blacksquare=\underline{R}\botimes\Z_\blacksquare[\underline{S}]$.
\end{enumerate}
\end{exmp}

Proposition 7.5 of \cite{condensed} shows that if $\A$ is an analytic ring, then we can define a subcategory of solid modules similar to how we defined solid abelian groups. This solid subcategory then satisfies properties analogous to those of $\SAb$. Examples of analytic rings include $\Z_\blacksquare$ (\cite[Theorem 5.8]{condensed}), $(A,\Z)_\blacksquare$ for $A$ discrete and $\Z_{p,\blacksquare}$ (\cite[Proposition 7.9]{condensed}). In fact, the proof for $\Z_{p,\blacksquare}$ being analytic extends easily to all profinite rings:

\begin{prop}
Let $R$ be a profinite ring. Then the pre-analytic ring $R_\blacksquare$ is analytic.
\begin{proof}
Exactly the same as that of \cite[Proposition 7.9]{condensed}. The key property of $R$ we need is that free solid $\underline{R}$-modules look like $\underline{R}\botimes$(free solid abelian groups); see (iii) from the example above.
\end{proof}
\end{prop}

As a consequence, \cite[Proposition 7.5]{condensed} tells us a lot about the category of solid $\underline{R}$-modules and its derived category. Recall that an $\underline{R}$-module $\M$ is solid if and only if $\Hom_{\underline{R}}(\underline{R}[\underline{S}]^\blacksquare,\M)=\Hom_{\underline{R}}(\underline{R}[\underline{S}],\M)$ for each $S\in\CHED$ (Corollary \ref{samesoliddef}). Similarly, let's call a complex $\C$ of $\underline{R}$-modules \emph{solid} if $\RHom_{\underline{R}}(\underline{R}[\underline{S}]^\blacksquare,\C)=\RHom_{\underline{R}}(\underline{R}[\underline{S}],\C)$ for each $S$. Then:

\begin{thm}
Let $R$ be a profinite ring.
\begin{enumerate}[label=(\roman*)]
\item The full subcategory $\mathbf{Solid}(\underline{R})$ of solid modules is an abelian subcategory of $\mathbf{CMod}(\underline{R})$ closed under all limits, colimits and extensions. The objects $\underline{R}[\underline{S}]^\blacksquare$ for $S\in\CHED$ form a family of compact projective generators. The inclusion $\mathbf{Solid}(\underline{R})\inj\mathbf{CMod}(\underline{R})$ admits a left adjoint $$\mathbf{CMod}(\underline{R})\to\mathbf{Solid}(\underline{R})\colon \M\mapsto\M^\blacksquare$$ which is the unique colimit-preserving extension of $\underline{R}[\underline{S}]\mapsto\underline{R}[\underline{S}]^\blacksquare$. If $R$ is commutative, there is a unique symmetric monoidal tensor product $\otimes^\blacksquare_{\underline{R}}$ on $\mathbf{Solid}(\underline{R})$ making the functor $\M\mapsto\M^\blacksquare$ symmetric monoidal.

\item The functor $D(\mathbf{Solid}(\underline{R}))\to D(\mathbf{CMod}(\underline{R}))$ is fully faithful, and its essential image is closed under all limits and colimits and given precisely by the solid complexes. If $\C$ is solid, then also $\RUHom_{\underline{R}}(\underline{R}[\underline{S}]^\blacksquare,\C)=\RUHom_{\underline{R}}(\underline{R}[\underline{S}],\C)$. A complex $\C\in D(\mathbf{CMod}(\underline{R}))$ lies in $D(\mathbf{Solid}(\underline{R}))$ if and only if all the $H^i(C)$ are solid. The inclusion $D(\mathbf{Solid}(\underline{R}))\inj D(\mathbf{CMod}(\underline{R}))$ admits a left adjoint $$D(\mathbf{CMod}(\underline{R}))\to D(\mathbf{Solid}(\underline{R}))\colon\C\mapsto\C^{L\blacksquare}$$ which is the left derived functor of $\M\mapsto\M^\blacksquare$. If $R$ is commutative, there is a unique symmetric monoidal tensor product $\otimes^{L\blacksquare}_{\underline{R}}$ on $D(\mathbf{Solid}(\underline{R}))$ making the functor $\C\mapsto\C^{L\blacksquare}$ symmetric monoidal. The functor $\otimes^{L\blacksquare}_{\underline{R}}$ is the left derived functor of $\otimes^\blacksquare_{\underline{R}}$.
\end{enumerate}
\begin{proof}
This is basically \cite[Proposition 7.5]{condensed} (except the very final statement), but we do have to check that the quoted proposition translates correctly to what we stated above. Also, note that most of (i) was stated previously as Proposition \ref{solidpr} and Corollary \ref{solidcor}. For (i), we know that the left adjoint to the inclusion (which would be written as $\M\mapsto\M\otimes_{\underline{R}}R_\blacksquare$ in \cite{condensed}) has to be $\M\mapsto\M^\blacksquare$, since the latter is a colimit-preserving extension of $\underline{R}[\underline{S}]\mapsto\underline{R}[\underline{S}]^\blacksquare$. Similarly, the solid tensor product (which would be written as $\otimes_{R_{\blacksquare}}$ in \cite{condensed}) has to be $(-\otimes^\blacksquare_{\underline{R}}-)=(-\otimes_{\underline{R}}-)^\blacksquare$ (see the paragraph before Proposition \ref{rtensors}).

As stated in \cite[Warning 7.6]{condensed}, it is not known whether the last statement of (ii) above is true in general (as far as the author is aware), but it is for profinite rings. To prove this, we need to show that for any $S,T\in\CHED$, $\underline{R}[\underline{S}]^\blacksquare\otimes^{L\blacksquare}_{\underline{R}}\underline{R}[\underline{T}]^\blacksquare=\underline{R}[\underline{S\times T}]^{L\blacksquare}$ is concentrated in degree 0. It is not obvious that $\underline{R}[\underline{S\times T}]^{L\blacksquare}=\underline{R}[\underline{S\times T}]^{\blacksquare}$ since $S\times T$ need not be extremally disconnected, but because we could have defined $\underline{R}[\underline{S}]^\blacksquare$ for all $S\in\Pro$ (rather than just for $S\in\CHED$), the argument of \cite[page 40]{analytic} gives the above isomorphism.
\end{proof}
\end{thm}

\appendix
\section{Basics of Condensed Modules}

This appendix serves as an introduction to condensed modules, especially for those who are somewhat new to sheaves and sites. No result below is original, but they are included to make the paper more self-contained. Readers who are familiar with sheaves and elementary condensed mathematics can safely ignore the appendix.

Let $\mathbf{TGrp}$, $\mathbf{TAb}$, $\mathbf{TRing}$ and $\mathbf{TMod}(R)$ denote respectively the categories of topological groups, topological abelian groups, topological rings, and topological modules over a topological ring $R$. Further, let $\mathbf{Pro}$, $\mathbf{PGrp}$, $\mathbf{PAb}$, $\mathbf{PRing}$ and $\mathbf{PMod}(R)$ denote respectively the categories of profinite spaces, profinite groups, profinite abelian groups, profinite rings, and profinite modules over a topological ring $R$. Recall the following from \cite[Definition 1.2]{condensed}:

\begin{defn}\label{condseta}
A \emph{condensed set/group/ring/\ldots} is a sheaf of sets/groups/rings/\ldots on the site of profinite spaces, with finite jointly surjective families of maps as covers. Equivalently, a condensed set/group/ring/\ldots is a functor $$\mathcal{T}\colon\Pro^{\mathrm{op}}\to\Set/\Grp/\Ring/\ldots$$ such that $\mathcal{T}(\varnothing)=*$ and:
\begin{enumerate}[label=(\roman*)]
\item For any profinite spaces $S_1, S_2$, the natural map $\mathcal{T}(S_1\sqcup S_2)\to\mathcal{T}(S_1)\times\mathcal{T}(S_2)$ is a bijection.
\item For any surjection $S'\sur S$ of profinite spaces with pullback $S'\times_SS'$ and its two projections $p_1, p_2$ to $S'$, the map $\mathcal{T}(S)\to\{x\in\mathcal{T}(S')\colon p_1^*(x)=p_2^*(x)\in\mathcal{T}(S'\times_SS')\}$ is a bijection.
\end{enumerate}

Given a condensed set/group/ring/\ldots $\mathcal{T}$, we call the collection of global sections $\mathcal{T}(*)$ its \emph{underlying set/group/ring/\ldots}.
\end{defn}

Note that for a condensed set $\mathcal{T}$, we only need the condition $\mathcal{T}(\varnothing)=*$ to ensure that $\mathcal{T}$ does not send every profinite space to the empty set. If we know that $\mathcal{T}(S)\neq\varnothing$ for some profinite space $S$, then (i) above implies that the map $\mathcal{T}(S)\to\mathcal{T}(S)\times\mathcal{T}(\varnothing)$ is a bijection, so we automatically get $\mathcal{T}(\varnothing)=*$. Thus, if $\mathcal{T}$ is a presheaf of groups or rings, for example, then (i) implies $\mathcal{T}(\varnothing)=*$.

There are some set-theoretic issues that are addressed in \cite{condensed} (Remarks 1.3, 1.4), which we will mostly ignore. Briefly, we should choose an uncountable strong limit cardinal $\kappa$ and use only profinite spaces of cardinality $<\kappa$ to define $\kappa$-small condensed sets. We then define the category of condensed sets as the (large) colimit of the categories of $\kappa$-small condensed sets. Let $\mathbf{CondSet}$, $\mathbf{CondGrp}$, $\mathbf{CondAb}$, $\mathbf{CondRing}$ and $\mathbf{CondMod}(R)$ denote respectively the categories of condensed sets, condensed groups, condensed abelian groups, condensed rings, and condensed modules over an abstract ring $R$ i.e.\ the codomain of $\mathcal{T}$ is $\ModR$.

By \cite[Proposition 2.7]{condensed}, the condensed categories above are equivalent to the categories of sheaves on the site of ($\kappa$-small) compact Hausdorff extremally disconnected spaces, with covers given by finite jointly surjective families of maps. Let $\CHED$ denote the category of compact Hausdorff extremally disconnected spaces and note that these spaces are also totally disconnected (so $\CHED$ is a full subcategory of $\Pro$). That is, we could have defined a condensed set/group/ring/\ldots as a functor $$\mathcal{T}\colon\CHED^\mathrm{op}\to\Set/\Grp/\Ring/\ldots$$ satisfying $\mathcal{T}(\varnothing)=*$ and the analogue of (i) in Definition \ref{condseta}. Note that the analogue of (ii) is automatic here.

Given a condensed ring $\R$, we will now define condensed $\R$-modules in the most natural way, as in \cite[Definition 2.3.22]{brink} (yet to appear).

\begin{defn}\label{condmoda}
Let $\R$ be a condensed ring. A \emph{(left) condensed $\R$-module} is a condensed abelian group $\M$ together with a natural transformation $\R\times\M\Rightarrow\M$ (viewed as condensed sets) such that for any $S\in\CHED$, the map $\R(S)\times\M(S)\to\M(S)$ makes $\M(S)\in\Ab$ a usual (left) $\R(S)$-module.

A \emph{morphism} between two condensed $\R$-modules is a natural transformation which is an $\R(S)$-module homomorphism for each $S$.
\end{defn}

Let $\CModRR$ denote the category of condensed $\R$-modules. (Note that we're not using the possibly better notation $\mathbf{CondMod}(\mathcal{R})$ to avoid confusion with the category $\CModR$ of condensed modules over an abstract ring $R$.)

\begin{lemma}\label{lemma1a}
Let $\R$ be a condensed ring. If $\M$ is a condensed $\R$-module then it is a condensed module over $\R(*)$ i.e.\ $\M(S)\in\mathbf{Mod}(\R(*))$ for any $S\in\CHED$. Moreover, this is functorial in $\M$, so we get a forgetful functor $\CModRR\to\mathbf{CondMod}(\R(*))$.
\begin{proof}
Given $S\in\CHED$, there is a unique map $S\to *$ and thus a ring homomorphism $\R(*)\to\R(S)$. So, we obtain a canonical functor $\mathbf{Mod}(\R(S))\to\mathbf{Mod}(\R(*))$. The rest follows.
\end{proof}
\end{lemma}

Given $T$ a topological space/group/ring/\ldots, there is a natural way to associate to it a condensed set/group/ring/\ldots (\cite[Example 1.5]{condensed}), namely we define the functor $\underline{T}$ via $\underline{T}(S)=C(S,T)$, the set of continuous maps from $S\in\CHED$ to $T$. The group/ring/\ldots structure on $C(S,T)$ is pointwise induced by that of $T$. This indeed makes $\underline{T}$ a sheaf.

\begin{lemma}\label{samesheafa}
Let $R$ be an abstract ring. The functor $M\mapsto\underline{M}$ from $\ModR$ to $\CModR$, where $M$ is equipped with the discrete topology, coincides with the constant sheaf functor (i.e.\ sheafification of the constant presheaf). In particular, $M\mapsto\underline{M}$ is left adjoint to the global sections functor $\mathcal{T}\mapsto\mathcal{T}(*)$.
\begin{proof}
Let $\Del(M)$ denote the constant presheaf with value $M$, and $\Del_0(M)$ denote the presheaf with $\Del_0(M)(S)=M$ for $S\neq\varnothing$ and $\Del_0(M)(\varnothing)=0$. The kernel of the obvious surjection of presheaves $\Del(M)\sur\Del_0(M)$ has sheafification equal 0. Since sheafification is exact, we conclude that the sheafification of $\Del_0(M)$ is the constant sheaf with value $M$ i.e.\ $\Del(M)^{\mathrm{sh}}$.

Given $M\in\ModR$ and $\varnothing\neq S\in\CHED$, there is an injection $M\inj C(S,M)$ sending $m\in M$ to the constant function with value $m$. This gives an injection of presheaves $\Del_0(M)\inj\underline{M}$ and hence an injection of sheaves $\Del(M)^{\mathrm{sh}}\inj\underline{M}$. To prove the lemma, it suffices to show that this map is surjective on each $S$. Take $f\in \underline{M}(S)=C(S,M)$, where $S\neq\varnothing$. Since $S\in\CHED$ is compact and $M$ is discrete, $f$ has finite image, say $\{m_1,\ldots,m_k\}\subseteq M$, with $f^{-1}(m_i)=S_i$. As $\Del(M)^{\mathrm{sh}}$ and $\underline{M}$ are both sheaves, we obtain a commutative diagram
\[\begin{tikzcd}
\Del(M)^{\mathrm{sh}}(S) \arrow[r] \arrow[d, "\cong"'] & \underline{M}(S) \arrow[d, "\cong"] \\
\prod\Del(M)^{\mathrm{sh}}(S_i) \arrow[r]              & \prod\underline{M}(S_i)            
\end{tikzcd}\]
Viewed as an element of $\prod\underline{M}(S_i)$, $f$ is just a constant function on each component $S_i$, so it certainly has a preimage in $\prod\Del(M)^{\mathrm{sh}}(S_i)$.
\end{proof}
\end{lemma}

Now let $R$ be a topological ring and $M$ a topological module over $R$, so we have $\underline{R}\in\CRing$ and $\underline{M}\in\CModR$, where the last $R$ is viewed as an abstract ring. We should expect $\underline{M}$ to be a condensed $\underline{R}$-module in the sense of Definition \ref{condmoda}, which is true:

\begin{lemma}\label{topmoda}
Let $R$ be a topological ring and $M$ a topological $R$-module. Then $\UM$ is naturally a condensed $\UR$-module. This is functorial, so we get an additive functor $\TModR\to\mathbf{CMod}(\underline{R})$. Moreover, the composition $\TModR\to\mathbf{CMod}(\underline{R})\to\mathbf{CondMod}(\underline{R}(*))=\CModR$, where the second map comes from Lemma \ref{lemma1a}, is just the functor $M\mapsto\underline{M}$.
\begin{proof}
Define a natural transformation $\phi\colon\UR\times\UM\Rightarrow\UM$ with component $\phi_S\colon C(S,R)\times C(S,M)\to C(S,M)$, for $S\in\CHED$, given by $\phi_S(f,g)(s)=f(s)g(s)$, where $f\in C(S,R), g\in C(S,M), s\in S$. Note that $\phi_S(f,g)$ is indeed continuous since it is the composition $S\overset{f\times g}{\longrightarrow} R\times M\to M$. The fact that this makes $C(S,M)$ an abstract $C(S,R)$-module follows easily from the fact that $M$ is an $R$-module.

It is also straightforward to show that $\phi$ is natural: given a map $\alpha\colon S'\to S$ in $\CHED$, the diagram
\[\begin{tikzcd}
{C(S,R)\times C(S,M)} \arrow[d] \arrow[r, "\phi_S"] & {C(S,M)} \arrow[d] \\
{C(S',R)\times C(S',M)} \arrow[r, "\phi_{S'}"']     & {C(S',M)}         
\end{tikzcd}\]
commutes as the pair $(f,g)\in C(S,R)\times C(S,M)$ is sent to $(s'\mapsto f[\alpha(s')]g[\alpha(s')])\in C(S',M)$ in either direction. Thus, $\UM$ is a condensed $\UR$-module. Functoriality is checked similarly.

Finally, if we view $\UM$ as a condensed module over $\UR(*)=R$ (with no topology) then the action of $R$ on $\UM(S)=C(S,M)$ is given by the pullback of the ring homomorphism $R\to C(S,R)$ which sends $r\in R$ to the constant map with value $r$. This is of course just the pointwise action of $R$ on $C(S,M)$.
\end{proof}
\end{lemma}

Remark: we could have just said that the condensed $\underline{R}$-module structure on $\underline{M}$ is obtained from the $R$-module structure on $M$ by applying the product-preserving functor $\underline{(-)}$ to all relevant diagrams, but we did the above explicitly as an illustration.

In \cite{condensed}, Theorem 2.2 states that the category of ($\kappa$-small) condensed abelian groups $\CAb$ is an abelian category satisfying Grothendieck's (AB3), (AB4), (AB5), (AB6), (AB3*) and (AB4*), just like the category $\Ab$ of abelian groups. Moreover, $\CAb$ is generated by compact projective objects. Recall that an object $M$ of an abelian category is \emph{compact} if $\Hom(M,-)$ commutes with filtered colimits. If $M$ is compact projective, then $\Hom(M,-)$ commutes with all colimits (as finite coproducts, filtered colimits and cokernels make all colimits in an abelian category). Since the category $\ModR$ of abstract $R$-modules also satisfies these Grothendieck's (AB) axioms and is very similar to $\Ab$, the category $\CModR$ of condensed modules over an abstract ring $R$ should behave like $\CAb$.

\begin{thm}\label{thm1a}
Let $R$ be an abstract ring. The category $\CModR$ of ($\kappa$-small) condensed modules over $R$ has all small limits and colimits and is an abelian category satisfying (AB3), (AB4), (AB5), (AB6), (AB3*) and (AB4*).
\begin{proof}
Exactly the same as in \cite[Theorem 2.2]{condensed}. The main point is that limits and colimits are computed pointwise here (limits are generally computed pointwise for abelian sheaves on a site, but colimits usually need to be further sheafified).
\end{proof}
\end{thm}

\begin{thm}\label{thm2a}
Let $R$ be an abstract ring. The category $\CModR$ of ($\kappa$-small) condensed modules over $R$ is generated by compact projective objects, that is, for any $\M\in\CModR$ there are some compact projective objects $\mathcal{P}_i\in\CModR$ and a surjection $\bigoplus\mathcal{P}_i\sur\M$. In particular, $\CModR$ has enough projectives.
\begin{proof}
Exactly the same as in \cite[Theorem 2.2]{condensed}. Note that the forgetful functor $\CModR\to\CSet$ has a left adjoint $\mathcal{T}\mapsto R[\mathcal{T}]$, where $R[\mathcal{T}]$ is the sheafification of the functor that sends $S\in\CHED$ to the abstract free $R$-module $R[\mathcal{T}(S)]$. The compact projective generators are then given by $R[\underline{S}]$ for $S\in\CHED$.
\end{proof}
\end{thm}

Remark: the category $\CAb$ of condensed abelian groups has no non-zero injectives, by \cite{condinj}.

\begin{exmp}\label{appeg}
\begin{enumerate}[label=(\roman*)]
\item Let $S_1,S_2\in\Pro$, so $\underline{S_i}\in\CSet$ and we can define $\Z[\underline{S_i}]\in\CAb$ as in the proof above. We then have that $\Z[\underline{S_1\sqcup S_2}]=\Z[\underline{S_1}]\times\Z[\underline{S_2}]$, which follows from the isomorphisms
\begin{eqnarray*}
\Hom_\CAb(\Z[\underline{S_1}]\times\Z[\underline{S_2}],\M)&=&\Hom_\CAb(\Z[\underline{S_1}]\oplus\Z[\underline{S_2}],\M)\\
&=&\Hom_\CAb(\Z[\underline{S_1}],\M)\times\Hom_\CAb(\Z[\underline{S_2}],\M)\\
&=&\Hom_{\CSet}(\underline{S_1},\M)\times\Hom_{\CSet}(\underline{S_2},\M)\\
&=&\M(S_1)\times\M(S_2)\\
&=&\M(S_1\sqcup S_2)\\
&=&\Hom_\CAb(\Z[\underline{S_1\sqcup S_2}],\M)
\end{eqnarray*}
(where $\M\in\CAb$) by using Yoneda.
\item Note that the category $\CSet$ of condensed sets has all (small) limits and colimits (see \cite[Lemma 4.33]{land}). By the exact same argument as (i), we have more generally, for $S_1,S_2\in\Pro$, that $\underline{S_1\sqcup S_2}=\underline{S_1}\coprod\underline{S_2}$, where the coproduct $\coprod$ is taken in condensed sets. This implies (i) because $\Z[-]\colon\CSet\to\CAb$, being a left adjoint, preserves colimits.
\item Part (i) in particular implies that if $S\in\Pro$ is finite, then $\Z[\underline{S}]=\bigoplus_S\underline{\Z}$, since $\Z[\underline{*}]=\Del(\Z)^\mathrm{sh}=\underline{\Z}$ (see Lemma \ref{samesheafa}).
\end{enumerate}
\end{exmp}

The category $\CAb$ of condensed abelian groups has a tensor product $\otimes$ (see \cite[page 13]{condensed}). Specifically, given $\M,\mathcal{N}\in\CAb$, $\M\otimes\mathcal{N}$ is defined as the sheafification of $S\mapsto\M(S)\otimes_{\Z}\mathcal{N}(S)$. Just like for abstract abelian groups, giving a morphism $\M\otimes\mathcal{N}\Rightarrow\mathcal{K}$ in $\CAb$ is equivalent to giving a bilinear map $\M\times\mathcal{N}\Rightarrow\mathcal{K}$ i.e.\ a natural transformation $\M\times\mathcal{N}\Rightarrow\mathcal{K}$ in $\CSet$ such that for each $S\in\CHED$, $\M(S)\times\mathcal{N}(S)\to\mathcal{K}(S)$ is bilinear in $\Ab$. This follows from the diagrams
\[\begin{tikzcd}
\mathcal{M}(S)\times\mathcal{N}(S) \arrow[d] \arrow[rr, two heads] \arrow[rdd] &                          & \mathcal{M}(S)\otimes\mathcal{N}(S) \arrow[d] \arrow[ldd] \\
\mathcal{M}(S')\times\mathcal{N}(S') \arrow[rr, two heads, shorten >=0.7ex] \arrow[rdd]         &                          & \mathcal{M}(S')\otimes\mathcal{N}(S') \arrow[ldd]         \\
                                                                               & \mathcal{K}(S) \arrow[d] &                                                           \\
                                                                               & \mathcal{K}(S')          &                                                          
\end{tikzcd}\]

\begin{lemma}\label{sheaftensora}
Let $\M,\mathcal{N}$ be presheaves of abelian groups on $\CHED$. Then sheafification commutes with tensor product, that is, $\M^\mathrm{sh}\otimes\mathcal{N}^\mathrm{sh}=(\M\otimes\mathcal{N})^\mathrm{sh}$, where the first $\otimes$ is in $\CAb$ and the second is for presheaves.
\begin{proof}
For any $\K\in\CAb$, we have the isomorphisms
\begin{eqnarray*}
\Hom_\CAb(\M^\mathrm{sh}\otimes\mathcal{N}^\mathrm{sh},\K)&=&\Bil(\M^\mathrm{sh}\times\mathcal{N}^\mathrm{sh},\K)\\
&=&\Bil(\M\times\mathcal{N},\K) \\
&=&\Hom_{[\CHED^{op},\Ab]}(\M\otimes\mathcal{N},\K)  \\
&=&\Hom_\CAb((\M\otimes\mathcal{N})^\mathrm{sh},\K).
\end{eqnarray*}
The only isomorphism which is not immediate is the second one. First note that sheafification commutes with the forgetful functor $U\colon\CAb\to\CSet$, so $U(\M^\mathrm{sh})$ is the sheafification of $\M$ viewed as a presheaf of sets. Since sheafification commutes with finite products, it is clear that $\Hom(\M^\mathrm{sh}\times\mathcal{N}^\mathrm{sh},\K)=\Hom(\M\times\mathcal{N},\K)$ in the category of presheaves/sheaves of sets. Lastly, notice that this isomorphism preserves bilinearity: the natural map $\M\times\mathcal{N}\Rightarrow\M^\mathrm{sh}\times\mathcal{N}^\mathrm{sh}$ is linear, so a bilinear map $\M^\mathrm{sh}\times\mathcal{N}^\mathrm{sh}\Rightarrow\K$ gives a bilinear map $\M\times\mathcal{N}\Rightarrow\K$. Conversely, bilinearity of a map $\M\times\mathcal{N}\Rightarrow\K$ can be expressed using some commutative diagrams, so sheafifying the diagrams tells us that the induced map $\M^\mathrm{sh}\times\mathcal{N}^\mathrm{sh}\Rightarrow\K$ is bilinear too.

The result follows by Yoneda.
\end{proof}
\end{lemma}

\begin{corr}
The tensor product $\otimes$ makes $\CAb$ a monoidal category with unit $\underline{\Z}$, the constant sheaf with value $\Z$.
\begin{proof}
Let us prove associativity of $\otimes$ as an illustration. By the previous lemma, $(\M\otimes\mathcal{N})\otimes\K$ can also be calculated as the sheafification of the presheaf $S\mapsto(\M(S)\otimes\mathcal{N}(S))\otimes\K(S)$. Now associativity immediately follows from the same fact in $\Ab$.
\end{proof}
\end{corr}

\begin{corr}
Let $\T_1,\T_2$ be condensed sets. Then $\Z[\T_1]\otimes\Z[\T_2]=\Z[\T_1\times\T_2]$ i.e.\ the functor $\Z[-]\colon(\CSet,\times)\to(\CAb,\otimes)$ is monoidal.
\begin{proof}
Both $\Z[\T_1]\otimes\Z[\T_2]$ and $\Z[\T_1\times\T_2]$ are sheafifications of $S\mapsto \Z[\T_1(S)]\otimes\Z[\T_2(S)]=\Z[(\T_1\times\T_2)(S)]$.
\end{proof}
\end{corr}

Using the above information, we could have also defined a condensed $\R$-module (Definition \ref{condmoda}) as follows:

\begin{prop}\label{condmod2a}
A condensed ring $\R$ is equivalently a monoid in $(\CAb,\otimes)$.

A condensed $\R$-module (as in Definition \ref{condmoda}) is equivalently an $\mathcal{R}$-module object in $\CAb$, that is, it is a condensed abelian group $\M$ together with a map $\R\otimes\M\Rightarrow\M$ in $\CAb$ which makes the following two diagrams commute:
\[\begin{tikzcd}
(\mathcal{R}\otimes\mathcal{R})\otimes\mathcal{M} \arrow[r, "\cong"] \arrow[d] & \mathcal{R}\otimes(\mathcal{R}\otimes\mathcal{M}) \arrow[r] & \mathcal{R}\otimes\mathcal{M} \arrow[d] &  & \underline{\mathbb{Z}}\otimes\mathcal{M} \arrow[r] \arrow[rd, "\cong"'] & \mathcal{R}\otimes\mathcal{M} \arrow[d] \\
\mathcal{R}\otimes\mathcal{M} \arrow[rr]                                       &                                                             & \mathcal{M}                             &  &                                                                         & \mathcal{M}                            
\end{tikzcd}\]

In the above language, a morphism $\M\Rightarrow\mathcal{N}$ in $\CModRR$ is a natural transformation of condensed abelian groups which makes the following diagram in $\CAb$ commute:
\[\begin{tikzcd}
\mathcal{R}\otimes\mathcal{M} \arrow[d] \arrow[r] & \mathcal{R}\otimes\mathcal{N} \arrow[d] \\
\mathcal{M} \arrow[r]                             & \mathcal{N}                            
\end{tikzcd}\]
\begin{proof}
We will only prove the last paragraph, as the rest is similar. For clarity, in this proof, let's write $\otimes$ for the tensor product of presheaves, and $\otimes^{\mathrm{sh}}$ for its sheafification. It is clear that a morphism $\M\Rightarrow\mathcal{N}$ in $\CModRR$ (as in Definition \ref{condmoda}) is a natural transformation which makes the top left square below commute:
\[\begin{tikzcd}
\mathcal{R}\otimes\mathcal{M} \arrow[r] \arrow[d] \arrow[rrd] & \mathcal{R}\otimes\mathcal{N} \arrow[d] \arrow[rrd] &                                                               &                                                     \\
\mathcal{M} \arrow[r] \arrow[rrd,equal]                             & \mathcal{N} \arrow[rrd,equal]                             & \mathcal{R}\otimes^{\mathrm{sh}}\mathcal{M} \arrow[d] \arrow[r] & \mathcal{R}\otimes^{\mathrm{sh}}\mathcal{N} \arrow[d] \\
                                                              &                                                     & \mathcal{M} \arrow[r,shorten >=1.2ex]                                         & \mathcal{N}                                        
\end{tikzcd}\]
If the bottom right square commutes, then so will the top left. The converse is also true by sheafifying the top left square.
\end{proof}
\end{prop}

Remark: our first definition of a condensed $\R$-module (Definition \ref{condmoda}) is also often called an ``$\R$-module object" in $(\CSet,\times)$. Whenever we use the term ``module object", we mean a module object over a monoid in a monoidal category. By the above, this should not cause any confusion.

\begin{lemma}\label{adj1a}
Let $\R$ be a condensed ring. The forgetful functor $\CModRR\to\CAb$ has a left adjoint $\R\otimes(-)$.
\begin{proof}
Given $\M\in\CAb$, we get a condensed abelian group $\R\otimes\M$. This is naturally an $\R$-module via the map $\R\otimes(\R\otimes\M)\cong(\R\otimes\R)\otimes\M\Rightarrow\R\otimes\M$. Let $\mathcal{N}\in\CModRR$. A map $\R\otimes\M\Rightarrow\mathcal{N}$ in $\CModRR$ will give a map $\M\cong\underline{\Z}\otimes\M\Rightarrow\R\otimes\M\Rightarrow\mathcal{N}$ in $\CAb$. Conversely, a map $\M\Rightarrow\mathcal{N}$ in $\CAb$ gives a map $\R\otimes\M\Rightarrow\R\otimes\mathcal{N}\Rightarrow\mathcal{N}$ which is easily checked to be a morphism in $\CModRR$.
\end{proof}
\end{lemma}

\begin{thm}\label{thm3a}
Let $\R$ be a condensed ring. The category $\CModRR$ of ($\kappa$-small) condensed $\R$-modules has all small limits and colimits and is an abelian category satisfying (AB3), (AB4), (AB5), (AB6), (AB3*) and (AB4*).
\begin{proof}
Given $i\mapsto\M_i$ with $\M_i\in\CModRR$, we try to define its limit $\mathcal{L}$ pointwise, so $\mathcal{L}(S)=\lim\M_i(S)$, a limit in $\mathbf{Mod}(\R(S))$. This makes $\mathcal{L}$ the limit of $\M_i$ in $\CAb$, by Theorem \ref{thm1a} and noting that the forgetful functor $\ModR\to\Ab$ preserves limits. There is then a unique map $\R\otimes\mathcal{L}\Rightarrow\mathcal{L}$ in $\CAb$ induced by the diagrams
\[\begin{tikzcd}
                                                    & \mathcal{R}\otimes\mathcal{L} \arrow[d, dashed] \arrow[ldd] \arrow[rdd] &                                          \\
                                                    & \mathcal{L} \arrow[ldd] \arrow[rdd]                                    &                                          \\
\mathcal{R}\otimes\mathcal{M}_i \arrow[rr] \arrow[d] &                                                                        & \mathcal{R}\otimes\mathcal{M}_j \arrow[d] \\
\mathcal{M}_i \arrow[rr]                            &                                                                        & \mathcal{M}_j                           
\end{tikzcd}\]
which makes $\mathcal{L}$ a condensed $\R$-module and hence the limit of $\M_i$ in $\CModRR$, by some diagram chasing. The point is that since $\mathcal{L}$ is a limit, we can show commutativity of diagrams such as those in Proposition \ref{condmod2a} by using the fact that there are unique maps into $\mathcal{L}$.

The proof for colimits is only slightly different. We try to define the colimit $\C$ of $i\mapsto\M_i$ pointwise, so that it is the colimit in $\CAb$. Since both sheafification and the tensor product $\otimes_{\Z}$ in $\Ab$ commute with colimits, we have $\R\otimes\C=\colim(\R\otimes\M_i)$ in $\CAb$. This induces a map $\R\otimes\C\Rightarrow\C$ as shown below which finishes the proof like before.
\[\begin{tikzcd}
\mathcal{R}\otimes\mathcal{M}_i \arrow[rr] \arrow[d] \arrow[rdd] &                                                 & \mathcal{R}\otimes\mathcal{M}_j \arrow[d] \arrow[ldd] \\
\mathcal{M}_i \arrow[rr] \arrow[rdd]                             &                                                 & \mathcal{M}_j \arrow[ldd]                             \\
                                                                 & \mathcal{R}\otimes\mathcal{C} \arrow[d, dashed] &                                                       \\
                                                                 & \mathcal{C}                                     &                                                      
\end{tikzcd}\]

\end{proof}
\end{thm}

\begin{thm}\label{thm4a}
Let $\R$ be a condensed ring. The category $\CModRR$ of ($\kappa$-small) condensed $\R$-modules is generated by compact projective objects. In particular, $\CModRR$ has enough projectives.
\begin{proof}
(See also \cite[Remarks 3.2.2]{brink}.) By Theorem \ref{thm2a}, $\CAb$ is generated by the compact projective objects $\Z[\underline{S}]$, $S\in\CHED$. Lemma \ref{adj1a} shows that the forgetful functor $U\colon\CModRR\to\CAb$ has a left adjoint $\R\otimes(-)$, while the proof of Theorem \ref{thm3a} shows that $U$ preserves all limits and colimits. For $\M\in\CModRR$, we have $\Hom(\R\otimes\Z[\underline{S}],\M)=\Hom(\Z[\underline{S}],U\M)$, so $\R\otimes\Z[\underline{S}]$ is compact projective in $\CModRR$ for each $S$. (Remark: $\R\otimes\Z[\underline{S}]$ should be called the free $\R$-module with basis $S$.)

Given $\M\in\CModRR$, we can forget to $\CAb$ and find a surjection $\bigoplus\Z[\underline{S_i}]\sur\M$. As $\R\otimes(-)$, being a left adoint, preserves all colimits, we have a surjection $\bigoplus(\R\otimes\Z[\underline{S_i}])\sur\R\otimes\M$ and hence a surjection $\bigoplus(\R\otimes\Z[\underline{S_i}])\sur\R\otimes\M\sur\M$.
\end{proof}
\end{thm}

Given a condensed ring $\R$, a right $\R$-module $\M$ and a left $\R$-module $\mathcal{N}$, there is a natural way to define their tensor product over $\R$, namely as the presheaf of abelian groups $(\M\otimes_{\R}\mathcal{N})(S)=\M(S)\otimes_{\R(S)}\mathcal{N}(S)$ for $S\in\CHED$. Note that this is actually a sheaf, because $$(\M(S_1)\times\M(S_2))\otimes_{\R(S_1)\times\R(S_2)}(\mathcal{N}(S_1)\times\mathcal{N}(S_2))=(\M(S_1)\otimes_{\R(S_1)}\mathcal{N}(S_1))\times(\M(S_2)\otimes_{\R(S_2)}\mathcal{N}(S_2)).$$ Explicitly, the isomorphism is given left-to-right by $(m_1,m_2)\otimes(n_1,n_2)\mapsto(m_1\otimes n_1,m_2\otimes n_2)$. Thus, the pointwise tensor product $\M\otimes_{\R}\mathcal{N}$ is a condensed abelian group.

Although we have taken $\R$, $\M$ and $\mathcal{N}$ above to be sheaves, there is a completely analogous (and standard) way to define $\R$-module objects for $\R$ a presheaf of rings on the site $\CHED$. We can then define $\M\otimes_{\R}\mathcal{N}$ just like above, which is now only a presheaf of abelian groups.

We have the following results analogous to Lemma \ref{sheaftensora} and the paragraph before it.

\begin{lemma}
Let $\R\in\CRing$. The tensor product $\M\otimes_{\R}\mathcal{N}$ represents $\R$-middle linear maps out of $\M\times\mathcal{N}$, that is, giving a map $\M\otimes_{\R}\mathcal{N}\Rightarrow\K$ in $\CAb$ is equivalent to giving a map $\M\times\mathcal{N}\Rightarrow\K$ in $\CSet$ such that for each $S$, $\M(S)\times\mathcal{N}(S)\to\K(S)$ is $\R(S)$-middle linear. A similar statement holds for presheaves instead of sheaves.

If $\R$ is a presheaf of rings on $\CHED$ and $\mathcal{N}$ is a (left) $\R$-module, then $\mathcal{N}^\mathrm{sh}$ is naturally a (left) $\R^\mathrm{sh}$-module. If also $\M$ is a right $\R$-module, then $\M^\mathrm{sh}\otimes_{\R^\mathrm{sh}}\mathcal{N}^\mathrm{sh}=(\M\otimes_{\R}\mathcal{N})^\mathrm{sh}$.
\begin{proof}
Similar to the proof of Lemma \ref{sheaftensora} and the paragraph before it. For the second paragraph, note that the $\R^\mathrm{sh}$-module structure on $\mathcal{N}^\mathrm{sh}$ is given by the sheafification of $\R\otimes\mathcal{N}\Rightarrow\mathcal{N}$.
\end{proof}
\end{lemma}

\begin{lemma}
The functor $R\mapsto\underline{R}$ from $\mathbf{Ring}$ to $\CRing$, where $R$ is equipped with the discrete topology, coincides with the constant sheaf functor. In particular, $R\mapsto\underline{R}$ is left adjoint to the global sections functor.
\begin{proof}
This is similar to Lemma \ref{samesheafa}. Recall that $\Del$ denotes the constant presheaf functor. By that lemma, the functor $R\mapsto\underline{R}$ agrees with the constant sheaf functor $R\mapsto\Del(R)^\mathrm{sh}$ when regarded as functors from $\Ab$ to $\CAb$ (hence when regarded as functors from $\mathbf{Set}$ to $\CSet$), so it suffices to show that their multiplicative structures also agree. This follows from the diagram below since both multiplications are induced from that of $R$.
\[\begin{tikzcd}
\Del(R)^{\mathrm{sh}}\times\Del(R)^{\mathrm{sh}} \arrow[r] \arrow[d, "\cong"'] & \Del(R)^{\mathrm{sh}} \arrow[d, "\cong"] \\
\underline{R}\times\underline{R} \arrow[r]                                     & \underline{R}                           
\end{tikzcd}\]
\end{proof}
\end{lemma}

\begin{corr}\label{cortensora}
Every condensed abelian group $\M$ has a unique $\underline{\Z}$-(bi)module structure. In particular, $\CAb=\mathbf{CMod}(\underline{\Z})$. Moreover, if $\M,\mathcal{N}\in\CAb$, then $\M\otimes\mathcal{N}=\M\otimes_{\underline{\Z}}\mathcal{N}$ (where the first $\otimes$ is $\otimes^\mathrm{sh}_{\Z}$ in $\CAb$, so a quick way to express the statement is `` $\otimes^\mathrm{sh}_{\Z}=\otimes_{\Z^\mathrm{sh}}$").
\begin{proof}
Again, let $\Del$ denote the constant presheaf functor. It is clear that every $\M\in\CAb$ has a unique $\Del(\Z)$-module structure, and that this induces some $\underline{\Z}$-module structure on $\M$ via sheafification (concretely, this action is given by the canonical isomorphism $\underline{\Z}\otimes\M\cong\M$). Now observe that every $\underline{\Z}$-module structure on $\M$ comes from the $\Del(\Z)$-module structure.

For the last part, notice that both functors are sheafifications of $S\mapsto \M(S)\otimes_{\Z}\mathcal{N}(S)=(\M\otimes_{\Del(\Z)}\mathcal{N})(S)$.
\end{proof}
\end{corr}

Remark: we could have proven Theorems \ref{thm1a} and \ref{thm2a} only for the case $R=\Z$ i.e.\ quote \cite[Theorem 2.2]{condensed}, since that is all we need to do everything above and to prove Theorems \ref{thm3a} and \ref{thm4a}. Note, by a similar argument to the above corollary, that $\CModR=\mathbf{CMod}(\underline{R})$, where $R$ is an abstract (discrete) ring. Then Theorems \ref{thm3a} and \ref{thm4a} immediately imply Theorems \ref{thm1a} and \ref{thm2a}. (In fact, it's easy to show that for $\R$ a presheaf of rings, the category of $\R$-modules which happen to be sheaves is canonically equivalent to the category of $\R^{\mathrm{sh}}$-modules.)

We used the term ``bimodule" in the statement of Corollary \ref{cortensora}, which hopefully has an obvious meaning. We will still spell it out below for clarity:

\begin{defn}
Let $\R,\R'\in\CRing$. A \emph{condensed $\R$-$\R'$-bimodule} is a condensed abelian group $\M$ which is both a left $\R$-module and a right $\R'$-module such that the following diagram commutes:
\[\begin{tikzcd}
\R\otimes \M\otimes \R' \arrow[r] \arrow[d] & \M\otimes\R' \arrow[d] \\
\R\otimes\M \arrow[r]                       & \M                    
\end{tikzcd}\]

A \emph{condensed $\R$-bimodule} is an $\R$-$\R$-bimodule as above.
\end{defn}

As explained on \cite[page 13]{condensed}, given $\M,\mathcal{N}\in\CAb$, the abelian group of morphisms $\Hom(\M,\mathcal{N})$ can be enriched to a condensed abelian group $\UHom(\M,\mathcal{N})$ by defining, for $S\in\CHED$, $\UHom(\M,\mathcal{N})(S)=\Hom_\CAb(\Z[\underline{S}]\otimes\M,\mathcal{N})\in\Ab$. This defines an internal Hom satisfying the Hom-tensor adjunction $\Hom_\CAb(\mathcal{K},\UHom(\M,\mathcal{N}))=\Hom_\CAb(\mathcal{K}\otimes\M,\mathcal{N})$ (for a proof, see \cite[Proposition 5.5]{land}). It is not necessarily obvious that the above formula for $\UHom(\M,\mathcal{N})(S)$ even defines a sheaf, but that follows from $\Z[\underline{S_1\sqcup S_2}]=\Z[\underline{S_1}]\times\Z[\underline{S_2}]$ (see Example \ref{appeg}(i)).

\begin{exmp}\label{uhomeg}
\begin{enumerate}[label=(\roman*)]
\item Given $\M,\mathcal{N}\in\CAb$, we can recover $\Hom_{\CAb}(\M,\mathcal{N})$ from $\UHom(\M,\mathcal{N})$ by taking global sections, since $\UHom(\M,\mathcal{N})(*)=\Hom(\Z[\underline{*}]\otimes\M,\mathcal{N})=\Hom(\M,\mathcal{N})$.
\item While the notation might be somewhat misleading, we do not have $\UHom(\M,\mathcal{N})=\underline{\Hom(\M,\Nc)}$ in general, where $\Hom(\M,\Nc)\in\Ab$ is given the discrete topology. To see this, note that $\UHom(\underline{\Z},\mathcal{N})=\Nc$, while $\underline{\Hom(\underline{\Z},\Nc)}=\underline{\Nc(*)}$. See however Proposition \ref{homunderline}.
\end{enumerate}
\end{exmp}

We might expect there to be a Hom-tensor adjunction for general condensed rings, similar to the abstract case. Let us establish this now. For $\R$ a condensed ring, we will write $\Hom_\R$ to mean $\Hom_{\CModRR}$, similar to the abstract case. When $\R=\underline{\Z}$, we shall often omit the ring and just write $\Hom$.

Let $\R,\R'$ be condensed rings, $\K$ a right $\R$-module, $\M$ an $\R$-$\R'$-bimodule, and $\mathcal{N}$ a right $\R'$-module. We wish to prove the adjunction $\Hom_{\R}(\mathcal{K},\UHom_{\R'}(\M,\mathcal{N}))=\Hom_{\R'}(\mathcal{K}\otimes_\R\M,\mathcal{N})$, so let us first prove the case $\R=\underline{\Z}$, where the adjunction reads $\Hom(\mathcal{K},\UHom_{\R'}(\M,\mathcal{N}))=\Hom_{\R'}(\mathcal{K}\otimes\M,\mathcal{N})$. Thus, we first need to define the condensed abelian group $\UHom_{\R'}(\M,\mathcal{N})$ and also find a right $\R'$-module structure on $\mathcal{K}\otimes\M$. The latter is easy, as we simply take the action to be $(\mathcal{K}\otimes\M)\otimes\R'\cong\mathcal{K}\otimes(\M\otimes\R')\to\mathcal{K}\otimes\M$. For the former, note that by taking $\K=\Z[\underline{S}]$, $S\in\CHED$ in the adjunction we wish to prove, we are forced to define $\UHom_{\R'}(\M,\mathcal{N})(S)=\Hom_{\R'}(\Z[\underline{S}]\otimes\M,\mathcal{N})\in\Ab$, which does define a sheaf. For a general condensed abelian group $\K$, we can express $\K$ as the cokernel of some map $\bigoplus_i\Z[\underline{S_i}]\to\bigoplus_j\Z[\underline{S_j}]$. Applying the functors $\Hom(-,\UHom_{\R'}(\M,\mathcal{N}))$ and $\Hom_{\R'}(-\otimes\M,\mathcal{N})$ (both of which turn colimits into limits) to this map proves the adjunction in the case $\R=\underline{\Z}$.

Now let $\R$ be a general condensed ring. We can then make the condensed abelian group $\UHom_{\R'}(\M,\mathcal{N})$ into a right $\R$-module by using the left action of $\R$ on $\M$; to be precise, the action of $\R$ on $\UHom_{\R'}(\M,\mathcal{N})$ is given by the image of the identity map $\UHom_{\R'}(\M,\mathcal{N})\to\UHom_{\R'}(\M,\mathcal{N})$ under the compositions
\begin{eqnarray*}
\Hom(\UHom_{\R'}(\M,\mathcal{N}),\UHom_{\R'}(\M,\mathcal{N}))&=&\Hom_{\R'}(\UHom_{\R'}(\M,\mathcal{N})\otimes\M,\mathcal{N}) \\
&\to&\Hom_{\R'}(\UHom_{\R'}(\M,\mathcal{N})\otimes\R\otimes\M,\mathcal{N}) \\
&=&\Hom(\UHom_{\R'}(\M,\mathcal{N})\otimes\R,\UHom_{\R'}(\M,\mathcal{N})).
\end{eqnarray*}
Arguing as in the previous paragraph proves the following:

\begin{prop}\label{homtensora}
Let $\R,\R'$ be condensed rings, $\K$ a right $\R$-module, $\M$ an $\R$-$\R'$-bimodule, and $\mathcal{N}$ a right $\R'$-module. Then we have the Hom-tensor adjunction $$\Hom_{\R}(\mathcal{K},\UHom_{\R'}(\M,\mathcal{N}))=\Hom_{\R'}(\mathcal{K}\otimes_\R\M,\mathcal{N}).$$
\end{prop}

\bibliographystyle{unsrt}

\end{document}